\pdfoutput=1
\RequirePackage{ifpdf}
\ifpdf % We~are running pdfTeX in pdf mode
\documentclass[pdftex]{sigma}
\else
\documentclass{sigma}
\fi

\usepackage{mathrsfs}
\usepackage{dsfont}
\usepackage{array}
\usepackage{tikz-cd}
\usetikzlibrary{shapes}
\usetikzlibrary[patterns]
\usetikzlibrary{decorations.pathreplacing}

\usepackage{eucal} %lettres calligraphiées

\newcommand{\Tor}{\mathrm{Tor}}

\newcommand{\vir}[1]{\big[#1\big]^\mathrm{vir}}
\newcommand{\fvir}[2]{[\![ #1 ]\!]^{#2}}
\newcommand{\ev}{\mathrm{ev}}
\newcommand{\pt}{\mathrm{pt}}

\newcommand{\m}[1]{\mfk\left[ #1 \right]}
\renewcommand{\d}[1]{\dfk\left[ #1 \right]}

\newcommand{\Pic}{\operatorname{Pic}}
\newcommand{\Alb}{\operatorname{Alb}}

\newcommand{\modulo}{ \operatorname{mod} }

\newcommand{\abcd}{\left(\begin{smallmatrix} a & b \\ c & d \\ \end{smallmatrix}\right)}

\newcommand{\PP}{\mathbb{P}}
\newcommand{\NN}{\mathbb{N}}
\newcommand{\ZZ}{\mathbb{Z}}
\newcommand{\RR}{\mathbb{R}}
\newcommand{\CC}{\mathbb{C}}

\newcommand{\QQ}{\mathbb{Q}}
\newcommand{\HH}{\mathbb{H}}

\renewcommand{\AA}{\mathbb{A}}

\newcommand{\btheta}{ {\underline{\theta}} }

\newcommand{\bfX}{ {\boldsymbol X} }
\newcommand{\bfD}{ {\boldsymbol D} }
\newcommand{\bfw}{ {\boldsymbol w} }
\newcommand{\bfm}{ {\boldsymbol m} }

\newcommand{\Dfk}{\mathfrak{D}}

\newcommand{\mfk}{\mathfrak{m}}

\newcommand{\dfk}{\mathfrak{d}}

\newcommand{\sfq}{\mathsf{q}}

\newcommand{\D}{\mathcal{D}}

\renewcommand{\L}{\mathcal{L}}
\newcommand{\N}{\mathcal{N}}

\newcommand{\M}{\mathcal{M}}

\newcommand{\X}{\mathcal{X}}

\newcommand{\V}{\mathcal{V}}
\newcommand{\C}{\mathcal{C}}

\renewcommand{\H}{\mathcal{H}}
\renewcommand{\O}{\mathcal{O}}

\newcommand{\Y}{\mathcal{Y}}

\renewcommand{\leq}{\leqslant}
\renewcommand{\geq}{\geqslant}

\renewcommand{\bar}{\overline}

\newcommand{\gen}[1]{\big\langle #1 \big\rangle}

\numberwithin{equation}{section}

\newtheorem{theo}{Theorem}[section]
\newtheorem*{theom}{Theorem}
\newtheorem{prop}[theo]{Proposition}
\newtheorem{coro}[theo]{Corollary}
\newtheorem{lem}[theo]{Lemma}

\theoremstyle{definition}
\newtheorem{defi}[theo]{Definition}

\newtheorem{remarkk}[theo]{Remark}
\newtheorem{expl}[theo]{Example}

\begin{document}

\allowdisplaybreaks

\newcommand{\arXivNumber}{2409.09472}

\renewcommand{\PaperNumber}{046}

\FirstPageHeading

\ShortArticleName{Correlated Gromov--Witten Invariants}

\ArticleName{Correlated Gromov--Witten Invariants}

\Author{Thomas BLOMME~$^{\rm a}$ and Francesca CAROCCI~$^{\rm b}$}

\AuthorNameForHeading{T.~Blomme and F.~Carocci}

\Address{$^{\rm a)}$~Universit\'e de Toulouse, Institut de Math\'ematiques, route de Toulouse, Toulouse 2000, France}
\EmailD{\href{mailto:thomas.blomme@math.univ-toulouse.fr}{thomas.blomme@math.univ-toulouse.fr}}

\Address{$^{\rm b)}$~Universit\`a di Roma Tor Vergata, Via della Ricerca Scientifica~1, Roma 00133, Italy}
\EmailD{\href{mailto:carocci@mat.uniroma2.it}{carocci@mat.uniroma2.it}}

%\ArticleDates{Received September 24, 2024, in final form June 09, 2025; Published online June 18, 2025}

\Abstract{We introduce a geometric refinement of Gromov--Witten invariants for $\mathbb P^1$-bundles relative to the natural fiberwise boundary structure. We call these refined invariant correlated Gromov--Witten invariants. Furthermore, we prove a refinement of the degeneration formula keeping track of the correlation. Finally, combining certain invariance properties of the correlated invariant, a local computation and the refined degeneration formula we follow floor diagram techniques to prove regularity results for the generating series of the invariants in the case of $\mathbb P^1$-bundles over elliptic curves. Such invariants are expected to play a role in the degeneration formula for reduced Gromov--Witten invariants for abelian and K3 surfaces.}


\Classification{14N35; 14N10; 14J26}

\section{Introduction}

\subsection{Setting}

In this paper, we introduce and study a refinement of the relative (or logarithmic) Gromov--Witten invariants \cite{abramovich2014stable,gross2013logarithmic,kim2010logarithmic,li2001stable} of the pair $(Y=\PP_X(\O\oplus L),D)$, where $D=D^+ + D^-=\PP_X(L)\oplus\PP_X(\O)$ and $X$ is a smooth projective variety over $\mathbb C$. We fix the following discrete data:
\begin{itemize}\itemsep=0pt
 \item $g$, $n$, $m$ are positive integers,
 \item $\bfw=(w_1,\dots,w_n)$ is a $n$-tuple of non-zero integers with $\sum w_i=0$; we moreover set \smash{$b(\bfw)=\sum_{w_i>0}w_i=-\sum_{w_i<0}w_i$};
 \item $\beta\in H_2(X,\ZZ)$ is an effective homology class.
\end{itemize}
Let $\M:=\M_{g,m}\bigl(Y|D^{\pm},\beta,\bfw\bigr)$ be the moduli space of logarithmic stable maps parametrizing
\[
f\colon\ \ (C,p_1,\dots,p_m,q_1,\dots,q_n)\to Y,
 \]
where $C$ is a stable curve of genus $g$ with image in class $\beta+b(\bfw)F\in H_2(Y,\ZZ)$, the marked point $q_i$ is mapped to $D^{\mathrm{sgn}(w_i)}$ with tangency order $|w_i|$, the marked points $p_i$ are mapped to the interior. The vector $\bfw$ is thus called \textit{tangency profile}.
The logarithmic (or equivalently relative~\cite{abramovich2014comparison}) Gromov--Witten invariants are defined integrating constraints pulled-back along evaluation maps~${
\ev \colon \M_{g,m}\bigl(Y|D^{\pm},\beta,\bfw\bigr) \longrightarrow X^n\times Y^m}$,
 against the virtual fundamental class~\smash{$\vir{\M}$}~\cite{gross2013logarithmic}.

 We propose a refinement for these GW invariants by remembering a discrete quantity depending on the position of the points mapped to the boundary divisor. This quantity is called a \textit{correlator}.

\subsection{The Albanese evaluation}
The idea of the refinement comes from reinterpreting the map to $\bigl(Y, D^\pm\bigr)=\bigl(\PP_X(\O\oplus L), D^\pm\bigr)$ as the data of a map $f\colon C\to X$ together with the choice of an isomorphism $f^*L\cong\O_C(\alpha)$ where $\alpha\in H^0\bigl(C,\bar{M}_C^{\rm gp}\bigr)$ is defined using the logarithmic structure on $C$. When $C$ is smooth, this condition takes the familiar form $f^*L\simeq\O_C(\sum w_i q_i)$.

\subsubsection[Case of X times P\^{}1]{Case of $\boldsymbol{ X\times\PP^1}$}

To simplify, let us momentarily assume that $L$ is the trivial bundle and that the source curve~$C$ is smooth.
 Then, reinterpreting $f\in\M$ as explained above, part of the data of a log stable map is the isomorphism
$
\O_C\bigl(\sum w_iq_i\bigr)\cong\O\in\Pic^0(C)$.
For a smooth curve, the Abel--Jacobi theorem ensures that $\Pic^0(C)\simeq\Alb(C)$, where $\Alb(C)$ is the Albanese variety of $C$.
 We denote by~${a_X\colon X\to\Alb(X)}$ the map from $X$ to its Albanese variety, unique up to translation. By the functoriality of the Albanese variety, $f$ induces a~morphism of Abelian varieties
$
 f_*\colon \Alb(C)\to\Alb(X)
$ compatible with the Albanese maps.

We denote by $a_\bfw$ the morphism
\[
 a_\bfw\colon\ (x_1,\dots,x_n)\in X^n\longmapsto \sum_{i=1}^n w_i a_X(x_i)\in\Alb(X).
 \]
As $\sum w_i=0$, the map $a_\bfw$ does not change if we compose $a_X$ with a translation. From $\O_C(\sum w_iq_i)\simeq\O\in\Pic^0(C)$, we get that
\begin{equation}\label{eq-f*f*}
 \sum_{i=1}^n w_i a_X(f(q_i)) =0\in\Alb(X).
\end{equation}
Equation \eqref{eq-f*f*} is a constraint on the relative position of the images $f(q_i)$. In other words, the~$X^n$ part of the evaluation map is not surjective, and truly has values in the subset $a_\bfw^{-1}(0)$.

\subsubsection{Case of a non-trivial bundle} If $L$ is not assumed to be trivial anymore, we prove in Lemma~\ref{lem:mapphibeta} that there is a natural map $\varphi_\beta\colon \Pic^0(X)\to\Alb(X)$ defined through Hodge theory, which only depend on the curve class. Then
Proposition \ref{prop-image-of-L-is-constant} ensures that equation \eqref{eq-f*f*} becomes the following:
\begin{equation}\label{eq-f*f*L}
\sum_{i=1}^n w_ia_X(f(q_i)) = \varphi_\beta(L)\in\Alb(X).
\end{equation}
In other words, the image depends only on the degree $\beta$ and the line bundle $L$ and not on $f$.

\subsection{Correlated classes and correlated GW invariants} We now assume that the tangency orders $w_i$ have a non-trivial common divisor $\delta$. We then consider the map
\[
a_{\bfw/\delta}\colon\ (x_1,\dots,x_n)\in X^n\longmapsto\sum_1^n\frac{w_i}{\delta}a_X(x_i)\in \Alb(X).
 \]
The above morphism induces a morphism from the space of log stable maps $\M$
\[\kappa^\delta(f)=\sum \frac{w_i}{\delta}a_X(f(q_i))\in\Alb(X).\]
It follows from equation \eqref{eq-f*f*L} that
$\delta\cdot\kappa^\delta(f)=\varphi_\beta(L)$.
Therefore, the map $\kappa^\delta$ takes values in the set of $\delta$-roots of $\varphi_\beta(L)$, denoted by $T^\delta(L,\beta)$. The latter is a torsor under the group of $\delta$-torsion elements in $\Alb(X)$, denoted by $\Tor_\delta(\Alb(X))$.
The value of $\kappa^\delta$ is called a \textit{correlator} and is denoted by $\theta$.
We want to refine the Gromov--Witten of $\bigl(Y,D^{\pm}\bigr)$ keeping track of the value of the correlator.

The discussion above assumes the source curve $C$ to be smooth, but the morphism $\kappa^\delta$
is also defined on the boundary of the moduli space $\M_{g,m}\bigl(Y|D^{\pm},\beta,\bfw\bigr)$, still with values in $T^\delta(L,\beta)$ (see Section~\ref{sec:correlationrefinement}).
This shows that
 $\M$ splits into distinct (possibly still disconnected) components~$\M^\theta$ indexed by the values of the correlator $\theta\in T^\delta(L,\beta)$; in particular the virtual class splits as a sum of the so-called \textit{correlated classes} \smash{$\vir{\M^\theta}$}. We define the \textit{full correlated class}~$\fvir{\M}{\delta}$
\[ \fvir{\M}{\delta}=\sum_{\theta\in T^\delta(L,\beta)} \vir{\M^\theta}\cdot(\theta)\in \QQ[\Alb(X)]\otimes H_\bullet(\M,\QQ), \]
which is an element of the group algebra with cycle coefficients with support on $T^\delta(L,\beta)\subset\Alb(X)$.

\textit{Correlated GW invariants} are obtained by capping these classes with some pullback by the evaluation map: if $\gamma\in H^\bullet(X^n\times Y^m,\QQ)$, we set
\[ \langle\langle\gamma\rangle\rangle^\delta_{g,\beta,\bfw} = \int_{\fvir{\M}{\delta}}\ev^*\gamma \in\QQ[\Alb(X)], \]
which is an element of the group algebra $\QQ[\Alb(X)]$ with support on the torsor $T^\delta(L,\beta)$. Distinct correlators may yield different Gromov--Witten invariants, as illustrated by the computations in the elliptic case.

It follows from the discussion above that the evaluation map truly has values in the set \smash{$a_\bfw^{-1}(\varphi_\beta(L))=\bigsqcup a_{\bfw/\delta}^{-1}(\theta)$}. A broader generalization of GW invariants would be obtained by allowing the pullback of cohomology classes from $a_\bfw^{-1}(\varphi_\beta(L))$ rather than $X^n$. Correlated invariants are a particular case of these invariants where we pull-back the classes $1_\theta\in H^0\bigl(a_\bfw^{-1}(\varphi_\beta(L)),\QQ\bigr)$, corresponding to the components $a_{\bfw/\delta}^{-1}(\theta)$ of $a_{\bfw}^{-1}(\varphi_\beta(L))$.

\subsection{Properties of the correlated classes}

\subsubsection{Deformation invariance}
 We prove that the correlated class is deformation invariant, so that we are right to speak about a refinement of the GW invariants. Here we are allowed to deform both the variety $X$ and the line bundle $L$, considering a family $\X\to S$ with a line bundle $\L\to\X$. Notice the correlators take values in the torsor inside the Albanese variety $\Alb_{\X/S}$ which is also deformed. Using the latter for some suitable choice of family, we furthermore prove the following.

\begin{theom}[Theorem~\ref{theo-torsor-invariance}]
 The correlated invariants \smash{$\langle\langle\gamma\rangle\rangle^\delta_{g,\beta,\bfw}$} are invariant under the action of the subgroup $\varphi_\beta\bigl(\Tor_\delta\bigl(\Pic^0(X)\bigr)\bigr)$.
\end{theom}

In the case where $X=E$ is an elliptic curve, we have $\Pic^0(E)\simeq\Alb(E)\simeq E$. If the chosen homology class is $\beta=a[E]$, the map $\varphi_\beta\colon E\to E$ is actually the multiplication by $a$, so that we have
\smash{$\varphi_{a[E]}(\Tor_\delta(E)) = \Tor_{\delta/\gcd(a,\delta)}(E)$}.
In the two extremal case, it means the following:
\begin{itemize}\itemsep=0pt
 \item If $a$ and $\delta$ are coprime, the subgroup is $\Tor_\delta(E)$, so that the classes $\vir{\M^\theta}$ all provide the same invariants.
 \item If $\delta|a$, the subgroup is trivial and we have no further information on the correlated invariants.
\end{itemize}
Specific to the elliptic case, choosing $L=\O$ so that $T^\delta(L,\beta)=\Tor_\delta(E)$, we furthermore prove that the invariant for the class \smash{$\vir{\M^\theta}$} only depends on $\theta$ through its order in $\Tor_\delta(E)$.

\subsubsection{Degeneration formula}
One of the most important tools in the computation of Gromov--Witten invariants is the degeneration formula \cite{kim2018degeneration,li2002degeneration,maulik2023logarithmic}. Our second contribution is a degeneration formula suited for the correlated classes.

We state two degeneration formulas. The first one is stated in Theorem~\ref{theo-refined-decomposition} and valid for any variety $X$. The drawback of the latter is that to be useful, one needs an expression for a~refinement of a~diagonal class. In the usual degeneration formula, such expression is provided by the~K\"unneth decomposition formula. In our situation, such a decomposition is not obvious to get.~The second formula deals with the particular case where the maps $H^\bullet(X^n,\QQ) \to  H^\bullet\bigl(a_{\bfw/\delta}^{-1}(\theta),\QQ\bigr)$ are surjective. Under this assumption, we recover an expression using the usual diagonal with its K\"unneth decomposition. To state the formula, we use the operator $\d{\frac{1}{\delta}}$ on $\QQ[\Alb(X)]$, which maps $(\theta)$ to the average of its $\delta$-roots.

\begin{theom}[Corollary \ref{theo-refined-decomposition-surjective-case}]
For any class $\gamma\in H^\bullet(X^n\times Y^m,\QQ)$, we have the following equality between intersection numbers:
\[ \ev^*\gamma\cap\fvir{\M_0}{\delta}
= \sum_\Gamma \frac{\prod w_e}{|\mathrm{Aut}(\Gamma)|}\bigl(\ev^*\gamma\cup\ev^*\bfD^\vee\bigr)\cap\d{\frac{1}{\delta/\delta_\Gamma}}\prod_V \fvir{\M_V}{\delta_\Gamma}, \]
where the sum is over degeneration graphs $($Definition {\rm\ref{defi-degeneration graph})}, $\bfD^\vee$ is the Poincar\'e dual of the class of the diagonal, and $\mathfrak{d}$ is the division operator in the group algebra.
\end{theom}

This technical surjectivity assumption is satisfied in the case where $X$ is an elliptic curve $E$ since in that case the map $a_{\bfw/\delta}^{-1}(\theta)\hookrightarrow X^n$ is just the inclusion of a subtorus.

\subsection{The elliptic case}
The remainder of the paper focuses on the easiest possible case, where $X$ is an elliptic curve. For simplicity in the introduction, we state the theorems for the case $L\cong\O$.

\subsubsection{Local invariants} In Section \ref{sec-computation-local-invariants}, we are able to provide an explicit computation of the following correlated Gromov--Witten invariants, called \textit{local GW invariants}
\[ \langle\langle \pt_0,1_{w_1},\pt_{w_2},\dots,\pt_{w_n} \rangle\rangle^\delta_{1,a[E],\bfw}. \]
For this specific choice of insertions, the invariant is a concrete count of curves which are in this situation covers of $E$. The enumeration is thus possible through a careful study of these covers.

\begin{theom}[Theorem~\ref{theo-expression-local-correlated-invariant}]
 The full local correlated invariant has the following expression:
 \[ \langle\langle\pt_0,1_{w_1},\pt_{w_2},\dots,\pt_{w_n}\rangle\rangle^\delta_{1,a[E],\bfw} = a^{n-1}w_1^2\boldsymbol{\sigma}_\delta(a), \]
 where $\boldsymbol{\sigma}_\delta(a)$ is a refinement of the sum of divisors function $\sigma(a)=\sum_{d|a}d$ $($see Section {\rm\ref{sec-computation-local-invariants})}.
\end{theom}

We provide an explicit expression for $\boldsymbol{\sigma}_\delta(a)$. These invariants are in a sense the most simple ones because they carry a unique interior point constraint. They enable the computation of invariants with more point constraints through the degeneration formula as seen in the next section.

\subsubsection{Floor diagrams and regularity}
Using the degeneration formula, the above local invariants are in fact the ground blocks enabling the computations of the more general correlated invariants $\langle\langle\pt^{n+g-1},1_{w_1},\dots,1_{w_n} \rangle\rangle^\delta_{g,a[E],\bfw}$. The degeneration formula reduces the computation to the enumeration of \textit{floor diagrams} with some new multiplicities.

Floor diagrams are a combinatorial tool introduced by Brugall\'e and Mikhalkin to~enumerate some planar curves subject to point conditions chosen in a stretched configuration \mbox{\cite{brugalle2007enumeration,brugalle2008floor}}. Many variations have since been introduced to deal with various situations, see \cite{blomme2022floor} for the version relevant for this paper.

Using the floor diagram algorithm, we are able to provide two regularity statements refining~\cite[Theorems 6.6 and 6.9]{blomme2022floor}. The first statement deals with the generating series in the class direction.

\begin{theom}[Theorem~\ref{theo-quasi-modularity}]
 Let $\bfw$ be a tangency profile of length $n$. The generating series
 \[ \sum_{a=1}^\infty \gen{\gen{\pt^{n+g-1},1_{w_1},\dots,1_{w_n}}}_{g,a[E],\bfw}^\delta \sfq^a \]
 is quasi-modular for the congruence group $\Gamma_0(\delta)$ with values in the group algebra $\CC[\Tor_\delta(E)]$.
\end{theom}

Compared to \cite[Theorem~6.9]{blomme2022floor}, the refinement loses regularity since quasi-modularity is only satisfied for the congruence subgroup $\Gamma_0(\delta)$ and not ${\rm SL}_2(\ZZ)$. The next regularity statement concerns the dependance in the tangency orders $\bfw$. Namely, we consider the following functions:
\[ N^\delta_{g,a}\colon\ \bfw\longmapsto \bigl\langle\bigl\langle \pt^{n+g-1},1_{w_1},\dots,1_{w_n} \bigl\rangle\bigr\rangle^\delta_{g,a[E],\bfw} \in\QQ[E], \]
defined on the set of tangency profile $\bfw$ divisible by $\delta$.

\begin{theom}[Theorem~\ref{theo-quasi-polynomiality}]
 For fixed $a$, $g$, $n$, $\delta$, the function $\bfw\mapsto N^\delta_{g,a}(\bfw)$ is piecewise polynomial.
\end{theom}

Theorem~\ref{theo-quasi-polynomiality} states that the piecewise polynomiality from \cite[Theorem 6.6]{blomme2022floor} survives the refinement. Such behavior has also been observed in a number of similar situations, such as double Hurwitz numbers and relative GW invariants of Hirzebruch surfaces \cite{ardila2017double}.

\subsection{Future directions}
Correlated GW invariants from the elliptic case are interesting due to their relation to the GW invariants of bielliptic surfaces, when caring about the torsion part in the homology group of a bielliptic surface. The GW invariants up to torsion have already be computed by the first author in \cite{blomme2024bielliptic}, and the use of correlated invariants is necessary for the torsion part.

Importantly, these invariants seem also to be necessary for the study of reduced GW invariants of abelian surfaces, and the mysterious multiple cover formula they should satisfy \cite{blomme2022abelian3,oberdieck2022gromov}. We plan to investigate this further in subsequent works. Indeed, looking at the reduced decomposition formula from \cite[Section 4.6]{maulik2010curves}, we notice the following. The hyperplane of values of the evaluation map is in fact disconnected when tangency orders along the gluing divisor have a common non-trivial divisor. The correlators give a way to concretely describe the connected components, only one of them containing the diagonal involved in the decomposition formula.

Many problems concerning these correlated invariants and their computation remain open, such as their computation in the presence of $\psi$-classes or $\lambda$-classes. It would also be interesting to see if we can have a correlated double-ramification cycle formula refining the formula from~\cite{janda2020double}.

\section{Logarithmic curves, line bundles and stable maps}
\subsection{Curves and line bundles}
We
refer the reader to \cite{ogus} for an extensive introduction to logarithmic structures and in particular for all the basic definitions which we do not recall in what follows.
\subsubsection{Logarithmic curves}
Let $(S,M_S)$ be a fine and saturated logarithmic scheme. A log curve over $S$ is a proper, integral and logarithmically smooth morphism $\pi\colon C\to S$ with connected, reduced one-dimensional (geometric) fibers.
Kato provided the following local characterization.
\begin{theo}[\cite{KatoF}]
For every geometric point $p\in C$ with $\pi(p)=s$, there exists an \'etale local neighbourhood of $p$ in $C$ with a strict \'etale morphism to
\begin{description}\itemsep=0pt
 \item[\emph{smooth point:}] $\mathbb A^1_S$ with log structure pulled back form the base, i.e., $\overline{M}_{C,p}=\overline M_{S,s}$;
 \item[\emph{marking:}] $\mathbb A^1_S$ with log structure generated by the $0_S$-section and $\pi^*M_S$, i.e., $\overline{M}_{C,p}=\overline M_{S,s}\oplus \mathbb N$;
 \item[\emph{node:}] $\operatorname{Spec}(\mathcal O_S[x,y]/xy-t)$ for some $t\in\mathcal O_S$, with log structure induced by the multiplication map $\mathbb A^2_S\to\mathbb A^1_S$ and $t\colon S\to \mathbb A^1$. In particular, $\overline{M}_{C,p}=\overline M_{S,s}\oplus \mathbb N e_x\oplus \mathbb N e_y\slash e_x+e_y=\delta$ for~${\delta\in \overline{M}_{S,s}}$ the so called \emph{smoothing parameter}.
\end{description}
\end{theo}

\subsubsection{Tropical curves}

To a family of logarithmic curves, we can associate a family of tropical curves.
Consider first $C\to (\operatorname{Spec}(k),Q)$ a logarithmic curve over a log point.

The \emph{tropicalization} $\Gamma$ is defined as the generalised cone complex obtained as the following colimit:
\[
\varinjlim_{\eta\rightsquigarrow\xi} \operatorname{Hom}\bigl(\overline{M}_{C,p},\mathbb R_{\geq 0}\bigr),
\]
where $\eta\rightsquigarrow\xi$ denote specialization maps inducing surjective morphisms $\overline{M}_{C,\xi}\to\overline{M}_{C,\eta}$. From the local description of the log structure on log smooth curves recalled above, it follows that $\Gamma$ comes equipped with a morphism of cone complexes
$\Gamma\to\sigma_Q $,
for $\sigma_Q=\operatorname{Hom}(Q,\RR_{\geqslant 0})$ the cone dual to the monoid $Q$.

Alternatively, $\Gamma$ can be thought as the dual graph of $C$ (having a vertex for each irreducible component, a leg for each marking and an edge $e\in\Gamma$ for each node) together with a metrization of the bounded edges with values in the base monoid $Q$. This means that we have an assignment~${l\colon E_b(\Gamma)\to Q}$ defined by associating to each bounded edge its smoothing parameter~${\delta_e\in Q}$.

\begin{remarkk}
Notice that for each map $Q\to\mathbb R_{\geq 0}$ we obtain a tropical curve in the classical sense. Since $Q\to\mathbb R_{\geq 0}$ correspond to a point of the cone $\sigma_Q$ dual to $Q$ (in the usual sense of toric geometry), we can think of a tropical curve metrized in $Q$ as a family of usual tropical curves on the dual cone.
\end{remarkk}
\begin{remarkk}\label{rem:tropfamilies}
The tropicalization makes sense over more general base log schemes (see also \cite[Definition~2.3.3.3]{{molchowiselog}}). Let $(C,M_C)\to(S,M_S)$ a log smooth curve;
for each $s\in S$, $\Gamma_s$ consists of the dual graph of $C_s$ metrized in $\overline{M}_{S,s}$.
 For each specialization $\eta\rightsquigarrow s$, we have a compatible diagram
\[
\begin{tikzcd}
\Gamma_s\ar[r,"c"]\ar[d,"\ell"] &\Gamma_{\eta}\ar[d,"\ell"]\\
\overline{M}_{S,s}\ar[r] & \overline{M}_{S,\eta},
\end{tikzcd}
\]
where $\overline{M}_{S,s}\to \overline{M}_{S,\eta}$ is a localization morphism followed by sharpification (quotient the invertibles) and $c$ is an edge contraction; more precisely, $c$ contracts the edges of $\Gamma_s$ whose length becomes zero in $ \overline{M}_{S,\eta}$.

\end{remarkk}
\subsubsection{Piecewise linear functions}
Unravelling the definition of global sections of the characteristic sheaf, we get a convenient description of $H^0\bigl(C,\overline{M}_C^{\rm gp}\bigr)$ in terms of continuous functions on the tropicalization $\Gamma$ which are linear with integral slope along the edges.
Let us first assume that $C$ is a log curve over a~${S=(\operatorname{Spec}(k),Q)}$.

Since $\overline{M}_C^{\rm gp}$ is a constructible sheaf, a
section is determined by giving $m_x\in\overline{M}_{C,x}^{\rm gp}$ at a generic point $x$ of each stratum,
such that the values $m_x$ are compatible with specialization, i.e., if $\xi$ generalize to $\eta$ then the image of $m_{\xi}$ under the surjection $\overline{M}_{C,\xi}^{\rm gp}\to \overline{M}_{C,\eta}^{\rm gp}$ is $m_{\eta}$.

Looking at Kato's classification, there are three type of strata on log smooth curves: the open strata corresponding to irreducible components of the curves excluding markings, i.e., vertices in the tropicalization; the closed strata corresponding to markings; the closed strata corresponding to nodes. Therefore, a section $\alpha\in H^0\bigl(C,\overline{M}_C^{\rm gp}\bigr)$ is the data of
\begin{itemize}\itemsep=0pt
 \item for each vertex $V$ of $\Gamma$ a value $\alpha(V)\in Q$;
 \item for each leg and edge an integer $w_e\in\mathbb Z$ (the slope along edge or leg in $\Gamma$);
 \item a compatibility condition for the values $\alpha(V)$, $\alpha(W)$ of vertices incident to the same edge~$e$: $\alpha(V)-\alpha(W)=\delta_e w_e$.
\end{itemize}
Indeed, suppose that $x_e$ is a node
of $C$ which generalizes to the generic points $\eta_V$ and $\eta_W$; then we have a map
$\overline{M}_{C,x_e}^{\rm gp}\to Q\times Q$
that is injective and gives an identification
\[\overline{M}_{C,x_e}^{\rm gp}=\{(a,b)\in Q\times Q \mid a-b\in\delta_e\mathbb Z\}.\]

We denote by $\operatorname{CL}(\Gamma,Q)$ the set of functions satisfying the above conditions, called \textit{cone-wise linear}.\footnote{We reserve the name piece-wise linear for those function which become linear along the edges after a subdivision.} The previous discussion can be summarized in the following equality (see also \cite[Section~3.3.3]{marcuswiselog} and references therein)
$H^0\bigl(C,\overline{M}_C^{\rm gp}\bigr)=\operatorname{CL}(\Gamma,Q)$.

 \begin{remarkk}
 Using the description of tropicalization for a family $(C,M_C)\xrightarrow{\pi}(S,M_S)$ of log curves given in Remark~\ref{rem:tropfamilies} and the above discussion, we can think of a section $ \alpha\in H^0\bigl(S,\pi_*\overline{M}_C^{\rm gp}\bigr)$ as a collection \smash{$\alpha_s\in \operatorname{CL}\bigl(\Gamma_s,\overline{M}_{S,s}^{\rm gp}\bigr)$} compatible under edge contraction, i.e., generalization.
 \end{remarkk}

\subsubsection{Line bundles}
Let $C\xrightarrow{\pi} S$ be a log smooth curve. We have a short exact sequence
\begin{equation}\label{eq:logseq}
0\to\mathcal O_C^{\times}\to M_C^{\rm gp}\to\overline{M}_C^{\rm gp}\to 0,\end{equation}
from which we obtain a long exact sequence of sheaves on $S$
\begin{align*}
 \cdots\to \pi_*M_C^{\rm gp}\to\pi_*\overline{M}_C^{\rm gp}\to
 R^1\pi_*\mathcal O_C^{\times}\to
 R^1\pi_*M_C^{\rm gp}\to R^1\pi_*\overline{M}_C^{\rm gp}\to\cdots.
\end{align*}
In particular, if the underlying scheme of $S$ is a point, so that $\pi_*\overline{M}^{\rm gp}_C=H^0\bigl(C,\overline{M}^{\rm gp}_C\bigr)$ and $R^1\pi_*\O_C^\times=H^1\bigl(C,\O_C^\times\bigr)=\Pic(C)$, we get a map
\smash{$
 H^0\bigl(C,\overline{M}_C^{\rm gp}\bigr)\to \Pic^0(C)$}, $
 \alpha \mapsto\mathcal O_C(-\alpha)$.

\begin{remarkk}
 The sign comes from observing that when $\alpha\in H^0\bigl(C,\overline{M}_C\bigr)$, the $\mathcal O_C^{\times}$-torsors of lifts of $\alpha$ to a section of $M_C$ comes equipped with a map to $\mathcal O_C$, coming from the structure morphism~${\epsilon \colon M_C\to\mathcal O_C}$ of the log structure. Up to changing $\alpha$ by $-\alpha$, we now deal with $\O_C(\alpha)$ instead.
\end{remarkk}

The restriction of the line bundle $\mathcal O_C(\alpha)$ to the irreducible components $C_V$ of $C$ can be explicitly described \cite[Proposition~3.3.3]{marcuswiselog}.
\begin{prop}
 Let $C\xrightarrow{\pi} S=(\operatorname{Spec}(k),Q)$ be a log curve and $\alpha\in\operatorname{CL}(\Gamma, Q)$, where $\Gamma$ denotes the tropicalization of $C$. Then for any vertex $V$,
 \[\mathcal O_C(\alpha)\rvert_{C_V}=\pi^*\mathcal O_S(\alpha(V))\otimes\mathcal O_{C_V}\biggl(\sum_{e\vdash V} s_e(\alpha) q_e\biggr),
 \]
where the sum is taken over all the edges and legs $e$ incident to $V$ and $s_e(\alpha)$ is the slope of~$\alpha$ along~$e$, oriented away from~$V$, where $q_e$ denotes the preimage of the node associated to the edge~$e$ on the considered component.
\end{prop}
\begin{remarkk}
 The result is more generally true if $C\to S$ is a logarithmic curve with constant degeneration.
\end{remarkk}

\subsection{Logarithmic line bundles and trivializations}
\begin{defi}[\cite{molchowiselog}]
 Let $C\to S$ be a log smooth curve. A \emph{logarithmic line bundle} $P$ on $C$ is a $M_C^{\rm gp}$-torsor such that its restriction $P\rvert_{C_s}$ to the geometric fibers has \emph{bounded monodromy}.
\end{defi}
We refer the reader to \cite{molchowiselog} for all explanations about the meaning and the necessity of the bounded monodromy. Here we only need that it is a condition on the $\overline{M}_C^{\rm gp}$-torsor $\overline{P}$ associated to $P$, and that it is automatically satisfied when $\overline{P}$ is isomorphic to the trivial $\overline{M}_C^{\rm gp}$-torsor.

In particular, we see from the long exact sequence
\[H^0\bigl(C,\overline{M}_C^{\rm gp}\bigr)\to\Pic(C)\to H^1\bigl(C,M_C^{\rm gp}\bigr)\to H^1\bigl(C,\overline{M}_C^{\rm gp}\bigr)\to\cdots\]
that for any $L\in \Pic(C)$ the associated $M_C^{\rm gp}$-torsor \smash{$L^{\times}\otimes_{\mathcal O_C^{\times}}M_C^{\rm gp}$} is in fact a logarithmic line bundle, where we denote by $L^{\times}$ the $\mathbb{G}_m$-torsor obtained from the line bundle removing its zero section.

\begin{defi}[{\cite[Definition~4.5.1]{marcuswiselog}}]
 Let $C\xrightarrow{\pi} S$ be a log smooth curve and $L$ a line bundle on $C$. A \emph{logarithmic trivialization} of $L$ is a section of \smash{$L^{\times}\otimes_{\mathcal O_C^{\times}}M_C^{\rm gp}$}.
\end{defi}
The definition is a modification of \cite[Definition~4.5.1]{marcuswiselog}, where the authors define logarithmic trivializations relative to the base. The difference comes from the fact that \cite{marcuswiselog} are interested in moduli of logarithmic maps to rubber targets, while we study the rigid case.
Similarly, the following discussion is parallel to
\cite[Proposition~4.5.3]{marcuswiselog}.

From \eqref{eq:logseq}, we see that the logarithmic line bundle \smash{$L^{\times}\otimes_{\mathcal O_C^{\times}}M_C^{\rm gp}$} can be thought as a \smash{$\mathcal O_C^{\times}$}-torsor over the trivial $\overline{M}_C^{\rm gp}$-torsor.
In particular, there is an exact sequence
\[0\to H^0\bigl(C,L^{\times}\otimes_{\mathcal O_C^{\times}}M_C^{\rm gp}\bigr)\to H^0\bigl(C, \overline{M}_C^{\rm gp}\bigr)\to\Pic(C),
 \]
where the second map is defined by
$\alpha\mapsto L(-\alpha)=L\otimes\O(-\alpha)$.
In other words, a logarithmic trivialization of $L$ is a trivialization of the line bundle $L(-\alpha)$ for some cone-wise linear function~$\alpha$.

\subsection[Logarithmic stable maps to P\^{}1-bundles and their degenerations]{Logarithmic stable maps to $\boldsymbol{\mathbb P^1}$-bundles and their degenerations}

This subsection follows closely \cite[Section~5]{marcuswiselog}, with two small differences. First, \cite{marcuswiselog} considers maps to rubber target, while we deal with the rigid case. Second, they study maps to $\mathbb{P}^1$ or to degeneration of $\mathbb P^1$ to a chain, while we consider maps to $\mathbb{P}^1$-bundles $Y=\mathbb P_X(\mathcal O_X\oplus L)$ and degenerations of the latter to certain chains $\Y_0=\bigl(Y_1,D_1^\pm\bigr)\cup \bigl(Y_2,D_2^\pm\bigr)\cup\dots\cup\bigl(Y_N,D_N^\pm\bigr)$ of~$\mathbb{P}^1$-bundles over $X$ obtained by gluing $Y_i$ to $Y_{i+1}$ along $D_i^+\cong D_{i+1}^-\cong X$.

For the benefit of the reader, we state the results in the form most convenient for us and sketch how to adapt the proofs of \cite[Section~5]{marcuswiselog} to our situation.

Let $X$ be a smooth projective variety and $L$ a line bundle on $X$.
For a fixed $\beta\in H_2(X,\mathbb Z)$ and a vector of integers $\bfw=(w_1,\dots,w_n)$ such that $\sum_{i=1}^n w_i=c_1(L)\cdot\beta$, we consider the moduli space
$\M := \M_{g,m}\bigl(Y|D^{\pm},\beta,\bfw\bigr)$
of logarithmic stable maps to $Y=\mathbb P_X(\mathcal O_X\oplus L)$ endowed with the divisorial log structure $D^++D^-=\mathbb P_X(\mathcal O_X) +\mathbb P_X(L)$ and with prescribed contact order $w_i$ at $q_i$ along $D^++D^-$ \cite{abramovich2014stable,gross2013logarithmic}.
Positive $w_i$ encode the tangencies with $D^+$ and the negative $w_i$ the tangencies with $D^-$.

This space compactifies the moduli space parametrizing maps from marked curves to $X$,
\[f\colon\ \bigl(C,p_1,\dots, p_m,q_1,\dots,q_n\bigr)\to X,\]
together with an isomorphism
\[\lambda\colon\ f^* L\biggl(-\sum a_i q_i\biggr)\cong\mathcal O_C.\]

Using the language of logarithmic trivializations introduced in the previous subsection, this description can be extended over the boundary.

\begin{defi}\label{def:logstablemaps}
 We define $\M(X,L)$ to be the stack over $\operatorname{LogSch}^{fs}$ whose objects are
 \begin{enumerate}\itemsep=0pt
 \item [(1)] A log smooth curve $C\to S$ with $n$ marked points with logarithmic structure and $m$ schematic marked points.
 \item [(2)] A stable map $f\colon (C,p_i,q_i)\to X$.
 \item [(3)] A logarithmic trivialization \smash{$\lambda\in H^0\bigl(C,f^*L^{\times}\otimes_{\mathcal O_C^{\times}}M_C^{\rm gp}\bigr)$}, such that the induced $\alpha\in H^0\bigl(C,\allowbreak\overline{M}_C^{\rm gp}\bigr)$ is locally on $C$ comparable with $0$ and the slope of $\alpha$ along the $i$-th marking is $w_i$.
 \end{enumerate}

We recall that a section \smash{$\alpha\in H^0\bigl(C,\overline{M}^{\rm gp}_C\bigr)$} is said comparable to $0$ if at each point $x\in C$ we have that either $\alpha(x)\geq 0$ or $\alpha(x)\leq 0$ in the partial order on $\overline{M}^{\rm gp}_{C,x}$ defined by the monoid $\overline{M}_{C,x}$. We refer the reader to \cite[Section~2.5]{RSPW2} for further explanations on the geometric meaning of this condition.
\end{defi}

\begin{prop}\label{prop:logstabelmaps}
The moduli space $\M(X,L)$ is equivalent as a stack over logarithmic schemes to the moduli space of logarithmic stable maps $\M$ of Gross--Siebert, Abramovich--Chen with target~${(Y=\PP_X(\mathcal O_X\oplus L),\mathcal M_Y)}$ and fixed tangencies $w_i$ along the markings $q_1,\dots,q_n$ {\rm\cite{abramovich2014stable,gross2013logarithmic}}. Here~$\mathcal M_Y$ is the divisorial log structure defined by the fiberwise $0$ and $\infty$ sections.
\end{prop}
\begin{proof}
The proof essentially consists in unravelling the definition of logarithmic stable map to a $\mathbb P^1$-bundle. For the case of log stable maps to toric varieties, the analogous result is proved, for example, in \cite[Proposition~2.5.1.1]{RSPW2}.

Let $\mathrm{Tot}(L)\xrightarrow{\pi} X$ be the total space of $L$ endowed with the zero section logarithmic structure~$M_L$. First, we argue that a logarithmic stable map
\smash{$C\xrightarrow{F} \mathrm{Tot}(L)$} is the same as a logarithmic trivialization $\lambda$ such that the associated $\alpha$ is a section of the log structure, i.e., $\alpha\in H^0\bigl(C,\bar{M}_C\bigr)$.
Indeed, by definition, a log map is a map $\underline{F}$ of the underlying schemes, namely $f\colon C\to X$ plus a section of $f^*L$, together with a morphism $\underline{F}^{-1}M_L\to M_C$ of the logarithmic structures. Using the Borne--Vistoli description of fine and saturated log structures \cite{BorneVistoli}, the latter is the same as a~morphism of the characteristic sheaves $\underline{F}^{-1}\bar{M}_L\to \bar{M}_C$ such that the maps to $\mathfrak{D}\mathrm{iv}_C$ commute. The first map comes from
$
 \bar{M}_L\cong\mathbb N_X\to \mathfrak{D}\mathrm{iv}_{\mathrm{Tot}(L)}$, $
1\to \mathcal O_L(X)\cong \pi^* L
$
pulling back along $\underline{F}$. It follows that a lift of the morphism of characteristic sheaves to a~morphism of the logarithmic structures correspond to an isomorphism $f^*L\cong^{\lambda}\mathcal O_C(\alpha)$ for $\alpha$ the section of $\bar{M}_C$ induced by the map on characteristic sheaves.

At this point we argue parallel to \cite[Section~2.5.1]{RSPW2}.

The $\mathbb P^1$-bundle $\mathbb P_X(\mathcal O_X\oplus L)$ with its fiberwise toric log structure can be constructed as the quotient of the rank 2 vector bundle $\mathcal O_X\oplus L$ minus the $0$ section for the fiberwise $\mathbb G_m$-action.

Then any logarithmic map $C\to \mathbb P_X(\mathcal O_X\oplus L)$ locally lifts to $\O_X\oplus L\setminus 0_X$.

By the discussion above, this lift can be represented by $a$, $b$ logarithmic trivializations of~$\mathcal O_C$ and $f^*L$ such that the associated $\bar{M}_C^{\rm gp}$ section $\bar{a}$ and $\bar{b}$ are in $\bar{M}_C$. Notice that the ratio~${\lambda=ab^{-1}}$ is a well defined $\mathbb G_m$-invariant section of the torsor \smash{$f^*L^{\times}\otimes_{\mathcal O_C^{\times}}M_C^{\rm gp}$}, namely a logarithmic trivialization of $f^*L$.

Notice however that not all the section $\alpha\in \rm{H}^0\bigl(C,\bar{M}_{C}^{\rm gp}\bigr)$ can arise this way: indeed, since $(a,b)$ must avoid the zero section, this means that locally around each point $x\in C$ either $\bar{a}_x\in \bar{M}_{C,x} $ or $\bar{b}_x\in\bar{M}_{C,x}$ is zero (i.e., lifts to an isomorphism of trivial $\mathbb G_m$-torsors). For $\alpha=\bar{a}-\bar{b}$, this translates precisely to being locally comparable to zero.
\end{proof}

\begin{remarkk}
If in Definition~\ref{def:logstablemaps}, we drop the request on the local comparability with zero of the section $\alpha\in H^0\bigl(C,\bar{M}_C^{\rm gp}\bigr)$, we would obtain the stack over $\operatorname{LogSch}^{fs}$ parametrizing log maps to the $\mathbb G_{log}$ bundle $\rm{Tot}(L)\times_{\mathbb G_m}\mathbb G_{log}$ associated with $L$. Similar moduli stack have been considered various times in the literature \cite{marcuswiselog,ranganathan2024logarithmic,RSPW2,ranganathan2020rational}.\footnote{We thank an anonymous referee whose comments suggested us to add this remark.}
\end{remarkk}

\begin{remarkk}
The definition and the characterization above can easily be adapted to the case of moduli of logarithmic stable maps to expansions \cite{kim2010logarithmic,li2001stable,maulik2023logarithmic}. Indeed, it is sufficient to add the condition that the values of $\alpha(V)$ are totally ordered; similar discussions appear in \cite{marcuswiselog,ranganathan2024logarithmic}.

In the classical language of \cite{li2001stable}, maps in the boundary of $\M$ which differ by the rubber action are identified. It is explained in detail in \cite{CN} how to see this from the logarithmic prospective.
\end{remarkk}

We are also interested in considering log maps to degenerations, i.e.,
\begin{itemize}\itemsep=0pt
\item Logarithmic maps to \smash{$\Y_t\xrightarrow{p}\mathbb A^1_t$} certain one parameter degenerations of $\mathbb P_{\X_{t\neq 0}}(\O\oplus\L)$ to
\[\Y_0=\bigl(Y_1,D_1^\pm\bigr)\cup\bigl(Y_2,D_2^\pm\bigr)\cup\dots\cup \bigl(Y_N,D_N^\pm\bigr),
\]
where $\bigl(Y_i,D^{\pm}_i\bigr)$ are $\mathbb P^1$-bundles over $\X_0$ glued along their infinity and zero sections;
$\mathbb A^1_t$ is endowed with its toric log structure and $\Y_t$ with the divisorial log structure defined by~${\Y_0+\mathcal{D}^+ +\mathcal{D}^-}$, where $\mathcal{D}^\pm$ are the $0$ and $\infty$-section of the family, and making $p$ log smooth.
\item Logarithmic maps to the log smooth scheme $\Y_0\to(\operatorname{Spec} k,\mathbb N)$ with the log structure pulled back from the family $\Y_t$
\end{itemize}

We recall the following definition.

\begin{defi}[{\cite[Definition~5.2.1]{marcuswiselog}}]
 A \emph{divided tropical line} over a logarithmic scheme $(S,M_S)$ is defined as follows. Let $P$ be an $M_S^{\rm gp}$-torsor and $\bar{P}$ the associated $\bar{M}_S^{\rm gp}$ torsors; fix $\gamma_0\leq\gamma_1\leq\cdots\leq\gamma_N$ a non empty collection of sections of $\bar{P}$ defined locally on $S$. Then $\bar{P}_{\gamma}\subset\bar{P}$ is defined to be the subfunctor whose local sections are comparable with all the $\gamma_i$.
\end{defi}
It is proved in \cite[Proposition~5.2.4]{marcuswiselog} that $P_{\gamma}:=\bar{P}_{\gamma}\times_{\bar{P}} P$ is a 2-marked family of semi-stable genus zero curve.

We then consider the following variation: fix $(\mathcal X, M_{\X})\xrightarrow{p} (S,M_S)$ with $X\to S$ smooth projective and $M_{\X}=p^*M_S$.
The two cases of interest listed above correspond to the case where~$(S,M_S)$ is $\mathbb A^1$ with its toric log structure or a standard log point.

Let $P\in\operatorname{LogPic}_{\X\slash S}$ and
$\bar{P}$ the induced $\bar{M}^{\rm gp}_{\mathcal{X}}$-torsor. Let $\gamma_0\leq\gamma_1\leq\cdots\leq\gamma_N$ be a non empty collection of sections of \smash{$\bar{P}$}. Notice that where the logarithmic structure of $\X$ is trivial $P$ is a line bundle in the usual sense and all the $\gamma_i$ coincide and are necessarily zero.
\begin{lem}\label{lem:tropicaline}
In the notation above, let $\bar{P}_{\gamma}\subseteq \bar{P}$ be the sub-funtcor of sections locally comparable with all the $\gamma_i$ and $P_{\gamma}=P\times_{\bar{P}}\bar{P}_{\gamma}$. Then $P_{\gamma}\to\X\to S$ is an expanded degeneration of $\Y=\mathbb P_{\X}(\L\oplus\O)$ for $\L$ the line bundle on $\X$ of lifts of $\gamma_0$ to a section of $P$.
\end{lem}
\begin{proof}
Since we assumed that the collection of sections is not empty, $\bar{P}$ admits a trivialization~${\gamma_0\colon\X\to\bar{P}}$ which we think as the zero of $\bar{M}_{\X}$ group. Notice that if the log structure on~$S$ and thus on $\X$ is generically trivial, then all the $\gamma_i$ coincide and are actually zero on the locus of $\X$ with trivial logarithmic structure. Then there exist a line bundle $\L_0$ representing the log line bundle $P$, namely the $\O_{\X}^\times$-torsor of lifts of $\gamma_0$ to a section of the log line bundle. We already proved in Proposition~\ref{prop:logstabelmaps} that the subfunctor of \smash{$P=L_0^\times\otimes_{\O_{\X}^\times}M_{\X}^{\rm gp}$} whose local sections are~comparable with $\gamma_0$ is the $\mathbb{P}^1$-bundle $\Y=\mathbb P_{\X}(\L_0\oplus\O)$ endowed with its boundary logarithmic structure.
Once a trivialization is fixed, we have that
\[
\gamma_i-\gamma_0=\delta_i\in H^0\bigl(\bar{M}_{\X}^{\rm gp}\bigr);
\] these are maps \smash{$\operatorname{Trop}(\X)=\operatorname{Trop}(S)\xrightarrow{\delta_i}\operatorname{Trop}(\Y)$}. Subdividing along the image, we obtained the desired degeneration.

 Alternatively, the statement can be proved following the steps of \cite[Proposition~5.2.4]{marcuswiselog}.
\end{proof}

Let $\mathcal Y_t\to\mathbb A^1$ a semi-stable degeneration coming form a subdivided tropical line; for example, this is the case if $\Y_t$ is constructed starting from $\mathbb P_{\X_t}(\O_{\X_t}\oplus\L_t)$ by successively blowing up the zero and the infinity section on the central fiber.
Denote by $\M\bigl(\Y_t\slash\mathbb A^1\bigr)$ the moduli space of logarithmic stable maps to $\Y_t\slash\mathbb A^1$ from families of log smooth curves with fixed genus $g$, number of markings $m$, contact order $(w_1,\dots, w_n)$ along $\D_t^-+\D_t^+$ and curve class $\beta$. Then we have the following.

\begin{coro}\label{prop:logmapsdegeneration}
 The moduli space $\M\bigl(\Y_t\slash\mathbb A^1\bigr)$ parametrizes the following data:
 \begin{enumerate}\itemsep=0pt
 \item[$(1)$] A diagram of logarithmic maps
 \[
 \begin{tikzcd}
 C\ar[r,"f"]\ar[d] & \X_t\ar[d]\\
 S\ar[r] &\mathbb A^1
 \end{tikzcd}
 \]
 for $C\to S$ a family of log smooth curves.
 \item[$(2)$] A section \smash{$\lambda\in H^0\bigl(C, f^* \mathcal L^{\times}\otimes_{\mathcal O^{\times}_{C}}M_{C}^{{\rm gp}}\bigr)$} such that the induces $\alpha\in H^0\bigl(C,\bar{M}_C^{\rm gp}\bigr)$ is locally comparable with $f^*\gamma_i$.
 \end{enumerate}
\end{coro}
\begin{remarkk}
 Parallel to before, if we want to instead consider logarithmic maps to expanded degenerations \cite{kim2010logarithmic,li2001stable,maulik2023logarithmic}, it suffices to further impose that for each geometric point $s$ the values of $\alpha(V)$ for $V\in\Gamma_s$ are totally ordered.
\end{remarkk}

\section{Correlated virtual class}

\subsection{Recollection on Albanese varieties}

\subsubsection{Complex setting}
For $X$ a proper, smooth variety over the complex numbers the Albanese variety was originally defined via transcendental methods using path integrals of closed holomorphic 1-forms. Let~$H_1(X,\mathbb Z)$ denote the torsion free part of the first homology group. There is an inclusion via the integration on paths:
\begin{align*}
 H_1(X,\mathbb Z)\longrightarrow H^0(X,\Omega_X)^*,\qquad
 \gamma\longmapsto \left(\omega\mapsto\int_{\gamma}\omega\right),
\end{align*}
where $ H^0(X,\Omega_X)^*$ denotes the dual to the vector space of holomorphic $1$-forms.
The Albanese variety of $X$ is the complex abelian variety obtained as the following quotient:
\[\Alb(X):=H^0(X,\Omega_X)^*\slash H_1(X,\mathbb Z). \]
We list below some results and properties of the Albanese variety and refer the reader to \cite[Section~11.11]{birkenhake} or \cite{mumford1974abelian} and references therein for all details and proofs.
\begin{enumerate}\itemsep=0pt
 \item[(1)] When $X$ is an Abelian variety, then $\Alb(X)\cong X$.
 \item[(2)] Let $x_0\in X$ a point, then there is an algebraic map, called the \emph{Albanese map}
 \begin{align*}
 a_{X}\colon \ X\longrightarrow \Alb(X),\qquad
 x\longmapsto\left(\omega\mapsto\int_{x_0}^x\omega\right) \mod H_1(X,\mathbb Z),
 \end{align*}
 sending $x_0$ to the identity element. This map is unique up to translation in the sense that choosing a different $x_0$ amounts to compose $a_X$ with a translation.
 \item[(3)] The Albanese map satisfies the following universal property:
 let $\varphi\colon X\to A$ be a map to an abelian variety, then there exists a unique homomorphism of abelian varieties $\widetilde{\varphi}\colon \Alb(X)\to A $ such that $\widetilde{\varphi}(0)=\varphi(x_0)$ and $\widetilde{\varphi}\circ a_{X}=\varphi$.
 \item[(4)] Let $X\xrightarrow{f} Y$ a morphism of algebraic varieties, then the map to the Albanese varieties are functorial, i.e., there exists a homomorphism of abelian varieties $f_*$ making the following diagram commute:
 \[
 \begin{tikzcd}
 X \arrow[r,"f"] \arrow[d,"a_X"']& Y \arrow[d,"a_Y"] \\
 \Alb(X) \arrow[r,"f_*"] & \Alb(Y),
 \end{tikzcd}
\]
 where $a_Y$ maps $f(x_0)$ to the identity.
 \item[(5)] The Albanese variety $\Alb(X)$ is the dual abelian variety of the Picard group
 \[\Pic^0(X):=H^1(X,\mathcal O_X)\slash H^1(X,\mathbb Z), \]
 i.e.,
 $\Alb(X)=\Pic^0\bigl(\Pic^0(X)\bigr)$.
 Identifying $\Pic^0(X)$ with the connected component of the identity of the group of line bundles, the functorial homomorphisms of Albanese variety~$f_*$ is the dual of the pullback map.
\end{enumerate}

\begin{remarkk}
 Notice that in the complex setting, $H^1(X,\O_X)$ can be seen as sheaf cohomology or alternatively as the Dolbeault cohomology group of $(0,1)$-forms.
\end{remarkk}

\subsubsection{Algebraic definition}
The interpretation of the Albanese variety as the dual abelian variety of \smash{$\operatorname{Pic}^0_{X\slash S}$}
allows a more algebraic definition.
 Indeed, when $\X\xrightarrow{g} S$ be a smooth, projective, morphism with connected geometric fibers over a scheme $S$, Grothendieck proved that there exists a smooth $S$-scheme, locally of finite presentation, $\Pic_{\X\slash S}$ representing the relative Picard functor and an abelian scheme (in particular, smooth and proper) $\Pic^0_{\X\slash S}$ over $S$ whose fibers $\Pic^0_{\X_s\slash k(s)}$ are the connected components of the identity in $\Pic_{\X_s\slash k(s)}$. See, for example, \cite[Section~9]{fantechi2006fundamental}, or \cite[Chapter~8]{bosch2012neron} and references therein.
Then by \cite[Theorem~3.3]{grothendieck1961technique} (see also \cite{milneAV,mochizuki2012topics}), we can consider the dual abelian scheme
\[\Alb^0_{\X\slash S}:=\Pic^0_{\Pic^0_{\X/S}\slash S}=\bigl(\Pic^0_{\X\slash S}\bigr)^\vee.\]
If $g\colon \X\to S$ has a section $\sigma_S\colon S\to\X$, then we get an Albanese map to $\Alb^0_{\X\slash S}$ constructed in the following way.

Let $\mathcal P$ be the \emph{unique} Poincar\'e line bundle on $\X\times\Pic^0_{\X\slash S}$ which trivializes along $\sigma_S$ in the sense of \cite[Section~4]{bosch2012neron}.

Then we define the Albanese morphism
\[
 a_{\X\slash S}\colon \ \X \longrightarrow \Alb^0_{\X\slash S},\qquad
 \big(U\xrightarrow{f} \X\big) \longmapsto (f\times\mathrm{id})^*\mathcal P.
 \]
Given \smash{$\X\xrightarrow{f} \Y$} a morphism of smooth, proper, geometrically connected $S$-varieties and $\sigma_S\colon S\allowbreak\to \X$ a section, it follows from the universal property of the Albanese scheme that there exists a~unique homomorphism $f_*$ of $S$-abelian schemes making the following diagram commute:
\[
\begin{tikzcd}
\X\ar[r,"f"]\ar[d,"a_{\X\slash S}"'] & \Y\ar[d,"a_{\Y\slash S}"]\\
\Alb^0_{\X\slash S}\ar[r,"f_*"] & \Alb^0_{\Y\slash S},
\end{tikzcd}
\]
where $a_{\X\slash S}$ and $a_{\Y\slash S}$ sends the sections $\sigma_S\colon S\to \X$ respectively $f\circ \sigma_S \colon S\to \Y$ to the identity element. Alternatively, $f_*$ can be described as the dual (in the category of abelian $S$-schemes, see \cite[Section~I.8]{milneAV}) of the pullback
$f^*\colon \Pic^0_{\Y\slash S}\to \Pic^0_{\X\slash S}$.

\subsubsection{Self-duality for smooth curves} As explained, for example, in \cite[Section~9.5.26]{fantechi2006fundamental}, when $\C\to S$ is a family of smooth, projective, geometrically connected genus $g>0$
curves over a Noetherian base, there exists an \emph{self-duality isomorphism}
$\Alb^0_{\C\slash S}\to \Pic^0_{\C\slash S}$.

 If we have a section $\sigma_0\colon S\to \C$, then the self-duality has a very explicit description using the Abel--Jacobi map $A_{\sigma_0}\colon\C\to\Pic^0_{\C/S}$ given by
\[
 A_{\mathcal \sigma_0}\colon \ \C \longrightarrow \Pic^0_{\C\slash S},\qquad
 (T\to\C) \longmapsto \mathcal O_{\C_T}(p_T -\sigma_0\rvert_T),
\]
where $p_T\colon T\to\C\times_S T$ is the induced section of $\C_T=\C\times_S T$.
The self-duality isomorphism is simply given by the pullback along the Abel--Jacobi map. The isomorphism does not depend on the choice, nor on the existence of the section, see, for example, \cite{esteves2002autoduality} for a proof.

Notice that $A_{\sigma_0}$ coincides with the Albanese map defined using the Poincar\'e line bundle on~${\X\times\Pic^0_{\X\slash S}}$ rigidified along the section $\sigma_0$.

\subsection{Refinement of the moduli spaces}

Let $\M=\M_{g,m}\bigl(Y|D^\pm,\beta,\bfw\bigr)$ be the moduli space of logarithmic stable maps to $Y=\mathbb P_X^1(\O\oplus L)$ considered in the previous section and let us now assume that $L\in\Pic^0(X)$: we have $c_1(L)=0$.

\subsubsection[Morphism to Alb(X)]{Morphism to $\boldsymbol{\operatorname{Alb}(X)}$}
Using the tangency orders $\bfw$, we consider the map
\[ a_\bfw\colon \ (x_i)\in X^n\longmapsto\sum_{i=1}^n w_i a_X(x_i)\in\Alb(X). \]
As $\sum w_i=0$, this map does not depend on the choice of a base point in $X$ and is uniquely defined. Indeed, for $z_0$, $z_0'$ two choices of base point and $a_X$, $a'_X$ the corresponding morphisms, we have that $a'_X= t_{a'_X(z_0)}\circ a_X$ and thus
\[
 \sum_{i=1}^n w_i a'_X(x_i)= \sum_{i=1}^n w_i (a_X(x_i) -a'_X(z_0))= \sum_{i=1}^n w_i a_X(x_i).
\]
Now, consider the morphism from the moduli space to the Albanese variety of $X$ by composing with the evaluation map:
\[ \kappa\colon\ \M\xrightarrow{\ev} X^n\xrightarrow{a_\bfw} \Alb(X)\]
defined by
\[%\label{eq:refinement morphism}
\bigl(\C\xrightarrow{f} X,(p_j)_{j=1}^m,(q_i)_{i=1}^n,\O_{\C}(\alpha)\cong^{\lambda} f^*L\bigr)\mapsto (x_i:=f\circ q_i)_{i=1}^n\mapsto \sum_{i=1}^n w_i a_X(x_i).
\]

Recall that $(q_i)_{i=1}^n$ are marking with logarithmic structure while $(p_j)_{j=1}^m$ are standard schematical points on $\C$.

\begin{lem}\label{lem:mapphibeta}
Let $\beta\in H_2(X,\ZZ)$ be the homology class of a complex curve. The bilinear map
\[ \phi \colon\ \alpha\otimes\omega\in H^{0,1}(X)\otimes H^{1,0}(X)\longmapsto \int_\beta \alpha\wedge \omega\in\CC \]
induces a morphism
\[ \varphi_\beta \colon\ \Pic^0(X)\longrightarrow \Alb(X). \]
\end{lem}
\begin{proof}
We use the isomorphisms $H^{0,1}(X)=H^1(X,\O_X)$ and $H^{1,0}(X)=H^0(X,\Omega_X)$ as well as the Hodge decomposition $H^2(X,\CC)=H^{0,1}(X)\oplus H^{1,0}(X)$. Through the above isomorphisms, the map $H^1(X,\ZZ)\to H^1(X,\O_X)$ yielding $\Pic^0(X)$ is actually the composition of $H^1(X,\ZZ)\to H^1(X,\CC)$ followed by the projection onto $H^{0,1}(X)$.
The bilinear map $\phi$ induces a~linear map~${\phi^*\colon H^{0,1}(X)\to H^{1,0}(X)^*}$. To prove that it descends to a map $\Pic^0(X)\to\Alb(X)$, we need to show that
\[ \phi^*\bigl(H^1(X,\ZZ)\bigr)\subset H_1(X,\ZZ)\subset H^0(X,\Omega_X)^*. \]
Let $\omega\in H^{1,0}(X)$ and $\gamma\in H^1(X,\ZZ)$. We write $\gamma\otimes\CC=\gamma^{0,1}+\gamma^{1,0}$ for its Hodge decomposition in $H^2(X,\CC)$. We have
\begin{align*}
\phi^*(\gamma)(\omega) = {}& \int_\beta \gamma^{0,1}\wedge\omega \qquad\text{(by definition)}\\
={} & \int_\beta \bigl(\gamma^{0,1}+\gamma^{1,0}\bigr)\wedge\omega
= \beta\cap (\gamma\cup\omega)
= (\beta\cap\gamma)\cap\omega,
\end{align*}
where the second equality follows from the fact that \smash{$\int_\beta \gamma^{1,0}\wedge\omega=0$}, since $\beta$ is the class of a~complex curve. The last line ensures that $\phi^*(\gamma)(\omega)$ is a period of $\omega$ since $\beta\cap\gamma\in H_1(X,\ZZ)$. Therefore, as desired, $H^1(X,\O_X)\to H^0(X,\Omega_X)^*$ descends to a map
\[ \varphi_\beta \colon\ H^1(X,\O_X)/H^1(X,\ZZ)\longrightarrow H^0(X,\Omega_X)^*/H_1(X,\ZZ). \tag*{\qed}
\]\renewcommand{\qed}{}
\end{proof}

\begin{remarkk}
Giving a non-degenerate positive bilinear map $H^1(X,\O_X)\otimes H^0(X,\Omega_X)\to\CC$ amounts to providing a polarization, i.e., an ample line bundle $\L$, on the Abelian variety $\Pic^0(X)$, see, for example,~\cite{mumford1974abelian}. In an Abelian variety, a polarization then induces a map
\[ A_\L\colon\ L\in\Pic^0(X)\longmapsto t_L^*\L\otimes \L^{-1}\in\Pic^0\bigl(\Pic^0(X)\bigr)=\Alb(X). \]
Assume $X$ is projective with hyperplane class $h\in H^2(X,\ZZ)$ giving a K\"ahler form on $X$. Then we have a polarization given by the non-degenerate positive bilinear form
\[ \phi\colon\ \alpha\otimes\omega\longmapsto \int_X h^{n-1}\wedge\alpha\wedge\omega. \]
Assume the class $\beta$ is obtained by intersecting sufficiently many hyperplane sections: $\beta$ is Poincar\'e dual to $h^{n-1}$. Then, the previously defined morphism~$\varphi_\beta$ is actually the morphism~$A_\L$ provided by the above polarization on $\Pic^0(X)$.
\end{remarkk}
We claim that the morphism $\kappa$ is in fact constant to $\varphi_\beta(L)$. To prove the latter, we start with the case of smooth curve, and then use it to deal with the case of nodal curves.

\begin{prop}\label{prop-image-of-L-is-constant}
The morphism $\kappa\colon\M\rightarrow \operatorname {Alb}^0(X)$ defined above is constant equal to $\varphi_\beta(L)$.
\end{prop}

\begin{proof}
Let us first look at a geometric point $s=\operatorname{Spec}(\mathbb C)\to\M$ of the moduli space such that~the source curve $C$ is smooth. We claim that in such case
$\kappa(s)=f_*f^*L$,
where~${\O_C \bigl(\sum_{i=1}^n w_i q_i\bigr)\cong f^*L}$ (by definition of the moduli functor) and
$f_*\colon\Pic^0(C)\to\operatorname{Alb}^0(X)$ is the morphism obtained dualizing $f^*\colon\Pic^0(X)\to\Pic^0(C)$ via the self-duality isomorphism of~$\Pic^0(C)$.

To see that, choose $p_0$ a point in $C$ and $z_0=f(p_0)$; then the functoriality of the Albanese morphism for compatible choices of base point gives
\[f_*(f^*L)=f_*\biggl(\O_C\biggl(\sum_{i=1}^n w_i q_i\biggr)\biggr)=f_*\biggl(\sum_{i=1}^n w_i A_{p_0}(q_i)\biggr)= \sum_{i=1}^n w_i a_X(x_i)\]
with the last equality following from the commutativity.

Since we are working over $\mathbb C$, we may use that for $X$ (resp.\ $C$), we have
 \[ \Pic^0(X)=H^1(X,\O_X)/H^1(X,\ZZ) \qquad \text{and}\qquad \Alb(X)=H^0(X,\Omega_X)^*/H_1(X,\ZZ). \]
Through Hodge theory, $H^1(X,\O_X)$ is isomorphic to the Dolbeault cohomology group $H^{0,1}(X)$. In particular, by functoriality, the pullback map $f^*\colon\Pic^0(X)\to\Pic^0(C)$ is in fact induced by~${f^*\colon H^1(X,\O_X)\to H^1(C,\O_C)}$, which is the pullback at the level of $(0,1)$-forms. Using that $H^0(X,\Omega_X)=H^{1,0}(X)$, the composition $f_* f^*$ is induced by the $\CC$-linear map
 \[ \varphi_C \colon\ H^{0,1}(X) \xrightarrow{f^*} H^{0,1}(C)\simeq H^{1,0}(C)^*\xrightarrow{f_*}H^{1,0}(X)^*, \]
where the second arrow is the dual to the pullback map for holomorphic forms. The middle isomorphism is provided by the Poincar\'e duality. Finally, if $\alpha\in H^{0,1}(X)$ is a $(0,1)$-form on $X$ and $\omega\in H^{1,0}(X)$ an holomorphic $1$-form, we have
 \[ \varphi_C(\alpha)(\omega) = \int_C f^*\alpha\wedge f^*\omega = \int_C f^*(\alpha\wedge\omega) = \int_{f(C)}\alpha\wedge\omega= \int_\beta \alpha\wedge\omega, \]
where the penultimate equality is push-pull formula.

Now, assume that $C$ is nodal; we write $\nu\colon\bigsqcup C_V\to C$ for the normalization map and $f_V\colon C_V\to X$ for the restriction of $f\circ\nu$ to the component $C_V$ of the normalization. Since each $C_V$ is now smooth, the discussion above tells us that the morphism
\[
f_{V\ast}f_V^*\colon\ \Pic^0(X)\to\Pic^0(C_V)\to\operatorname{Alb}^0(X)
\]
is equal to $\varphi_{\beta_V}$, where $\beta_V=f_{V\ast}[C_V]$. Furthermore, since $\sum \beta_V=\beta$, we have
\begin{align*}
\varphi_{\beta}=\sum \varphi_{\beta_V}\colon\ \Pic^0(X) \to\prod_V\Pic^0(C_V)\to\operatorname{Alb}^0(X),\qquad
 L \to \nu^*f^*L\to \sum_V f_{V\ast}(f_V^*L).
\end{align*}

To conclude the proof, we only need to argue that
\[ \sum_V f_{V\ast}(f_V^*L)=\kappa(s).\]

When $C$ is singular, $f^*L\cong^{\lambda}\O_C(\alpha)$ for $\alpha\in H^0\bigl(C, \bar{M}_C^{\rm gp}\bigr)$ such that
\[\O_C(\alpha)\rvert_{C_V}=\O_{C_V}\biggl(\sum_{e\vdash V}s_{e,V}(\alpha)q_e\biggr),\]
satisfying the balancing condition $\sum_{e\vdash V}s_{e,V}(\alpha)=0$. Here $s_{e,V}(\alpha)$ is the outgoing slope of $\alpha$ at~$V$. Notice that $\alpha$ induces an orientation on the edges of $\Gamma$. We denote with $w_e$ the slope along an edge with this induced orientation. Then if $V$, $W$ are the two vertices adjacent to a~bounded edge $e$, we will have $s_{e,V}(\alpha)=-s_{e,W}(\alpha)$ with one of the two being $w_e$.

 We know from the analysis for the smooth case that
 \[f_{V\ast}(f_V^*L)=\sum_{e\vdash V}s_{e,V}(\alpha) a_X(f_V\circ q_i).\]

 From the observation above and the fact that $f_V(q_e)=f_W(q_e)$ since the maps are induced by the maps on the nodal curves, it follows that when we consider $\sum_Vf_{V\ast}(f_V^*L)$ the contributions coming from the nodes cancel out and we obtain the required identity.
\end{proof}

\begin{expl}
 Assume $X=E$ is an elliptic curve. We have $\Alb(E)\cong\Pic^0(E)\cong E$. For the homology class $\beta=a[E]$, the associated map $\varphi_\beta$ is just the multiplication by $a$.
\end{expl}

\begin{remarkk}
For a family \smash{$C\slash S\xrightarrow{F} X\times S$} with non-smooth family of source curves, it is not obvious how to define
\[F_*\colon\ \Pic^{[0]}_{C\slash S}\to\operatorname{Alb}^0(X)\times S.\]
Indeed, \smash{$\Pic^{[0]}_{C\slash S}\to S$} is not an abelian scheme, as it is no longer proper, and we do not have an Abel--Jacobi map. There are two possible ways to fix the situations.
\begin{itemize}\itemsep=0pt
\item We can consider the compactified Jacobian \smash{$\Pic^{[0]}_{C\slash S}\subseteq \overline{\mathbb J}_{C/S}^{\mathcal E_{\pi},\sigma}$} parametrizing rank one, torsion free simple $(\mathcal E_{\pi},\sigma)$-quasi-stable sheaves on the fibers of $C\to S$ for $\mathcal E_{\pi}$ the canonical degree~$0$ polarization on $C\slash S$ in the sense of \cite[Section~2.1]{melo2015compactifications}. We refer the reader to~\cite{melo2015compactifications} and to~\cite{MeloRapagnettaViviani} for existence and properties of the compactifications.
Then a family version of the argument in \cite[Section~5]{MeloRapagnettaViviani} shows that there is an isomorphism of $S$-group
\[
\Pic^{[0]}_{C\slash S}\to \Pic^0\bigl(\overline{\mathbb J}_{C/S}^{\mathcal E_{\pi},\sigma}\bigr)
\]
and we can then define $F_*$ as the dual of $F^*$ composed with the natural inclusion of \smash{$\Pic^{[0]}_{C\slash S}$} with the compactified Jacobian.

\item Alternatively, one should in future be able to use the theory of logarithmic abelian varieties and of logarithmic Picard group developed by Molcho--Wise \cite{molchowiselog}. They show that~$\operatorname{LogPic}_{C^v\slash S}\to S$ (where the superscript $C^v$ means we are considering $C\to S$ endowed with the vertical log structure)
is a proper group stack over $\operatorname{LogSch}$
whose fibers are logarithmic abelian varieties in the sense of Kajiwara, Kato and Nakayama \cite{kajiwara2008logarithmic}.
It is expected that many classical results about abelian varieties admit generalization for logarithmic abelian varieties. In this spirit, in a forthcoming paper \cite{delignepair} Molcho--Ulirsch--Wise will show that $\operatorname{LogPic}_{C^v\slash S}$ is self dual in the category of logarithmic Abelian variety, admits a~Abel--Jacobi map and is universal for morphism to abelian groups over log schemes.
Once such a theory will be fully developed, considering the composition
 \[
 \Pic^{[0]}_{C\slash S}\hookrightarrow \operatorname{LogPic}_{C^v\slash S}\xrightarrow{F_*}\Alb(X)\times S
 \]
we would obtain the desired extension.
\end{itemize}

Once the extension $F_*$ is defined, one can look at the closed points to give an explicit description of $F_*F^*L$ and verify that this coincide with $\sum_{i=1}^n w_i a_X(f\circ q_i)$.
Since we do not need the extension at the level of universal Jacobian but only the map from $\M$, we do not fill in the technical details left out in this remark.
\end{remarkk}

\subsubsection{Refinement by correlation}\label{sec:correlationrefinement}
%\tnote{Still two moduli spaces ...}
Let us go back to the setting of $\M(X,L)$ where we fix the ramification profile. To simplify the notation, for a given point in $\M$, we set $x_i=f\circ q_i$. We are interested in the case where the contact orders have a non-trivial common divisor, i.e., there is some $1\neq\delta|\gcd(w_i)$; the we can consider the morphism
\[ a_{\bfw/\delta}\colon\ (x_i)\in X^n\longmapsto \sum_{i=1}^n \frac{w_i}{\delta}a_X(x_i). \]
Composing with the evaluation map, we define
\begin{gather}
 \kappa^\delta \colon\ \M \longrightarrow T^\delta(L,\beta)\subset\Alb(X), \nonumber\\
 \bigl(f\colon C\to X,(p_j)_{j=1}^m,(q_i)_{i=1}^n, \O_C(\alpha)\rvert_{[f]}\cong^{\lambda}f^*L\bigr) \longmapsto \sum_{i=1}^n\frac{w_i}{\delta} a_X (x_i)
.\label{eq:maptotorsion}
\end{gather}
In fact, $\kappa^\delta$ composed with the multiplication by $\delta$ is the previously defined $\kappa$. As $\kappa$ is constant equal to $\varphi_\beta(L)$ by Proposition \ref{prop-image-of-L-is-constant}, this ensures that $\kappa^\delta$ has values in $ T^\delta(L,\beta)$, the set of $\delta$-roots of $\varphi_\beta(L)$. The latter is a torsor under the finite group $\Tor_\delta(\Alb(X))$ of $\delta$-torsion elements in~$\Alb(X)$. In the particular case where $L=\mathcal O_X$, $T^\delta(X,\O_X)$ is exactly the set of $\delta$-torsion elements in the Albanese variety $\Alb(X)$.

As $L$ does not possess any canonical $\delta$-root, Proposition \ref{prop-image-of-L-is-constant} no longer applies for $\kappa^\delta$, which has no reason to be constant anymore. The morphism $\kappa^\delta$ defined in \eqref{eq:maptotorsion} thus determines a decomposition of the moduli space of logarithmic maps to $\mathbb P_X(\mathcal O_X\oplus L)$ into components, according to the value taken by $\kappa^\delta$. This refinement is called \textit{correlation}, as it stems from relations between the images of the points inside the Albanese variety. The image $\kappa^\delta(f)$ is called a \textit{correlator} and is usually denoted by $\theta$.

\begin{remarkk}
 The preimages $\bigl(\kappa^\delta\bigr)^{-1}(\theta)$ of elements $\theta\in T^\delta(L,\beta)$ are called components of the moduli space but they are not necessarily connected.
\end{remarkk}

In particular, it follows from the discussion above that we have the following.
\begin{prop}
 Let $\M$ be the moduli space of logarithmic stable maps to $Y=\PP_X(\mathcal O\oplus L)$ with the divisorial log structure defined by $D^- + D^+$. Assume that $L\in\Pic^0(X)$ and that $\delta|\gcd(w_i)$.
 Then we have a splitting of the virtual class
 \[ \vir{\M}=\sum_{\theta\in T^\delta(L,\beta)}\vir{\M^\theta},\]
 where $T^\delta(L,\beta)$ is the $\operatorname{Tor}_{\delta}(\Alb(X))$-torsor of $\delta$-roots of $\varphi_\beta(L)$.
\end{prop}

\begin{defi}\label{def:fullrefined}
We encompass the information of the correlated virtual classes in the \textit{full correlated virtual class}
\[ \fvir{\M}{\delta} = \sum_{\theta\in T^\delta(L,\beta)} \vir{\M^\theta}\cdot (\theta),\]
which is an element in the group algebra $\QQ[\Alb(X)]$ with cycle coefficients and support over the torsor $T^\delta(L,\beta)$.
\end{defi}

\subsubsection{Correlated Gromov--Witten invariants}
The projection $\PP_X(\O\oplus L)\to X$ identifies $D^-$ and $D^+$ with $X$. The evaluation map hence takes the following form:
\[
\ev\colon\ \M_{g,m}\bigl(Y|D^\pm,\beta,\bfw\bigr) \longrightarrow Y^m\times X^n,
 \]
we integrate pullback of cohomology classes over the full correlated virtual class to get \textit{correlated Gromov--Witten invariants}: for $\gamma_1,\dots,\gamma_m\in H^\bullet(Y,\QQ)$ and $\widetilde{\gamma}_1,\dots,\widetilde{\gamma}_n\in H^\bullet(X,\QQ)$, we set
\[ \bigl\langle\bigl\langle\gamma_1,\dots,\gamma_m,\widetilde{\gamma}_1,\dots,\widetilde{\gamma}_n\bigr\rangle\bigr\rangle^\delta_{g,\beta,\bfw} = \int_{\fvir{\M}{\delta}} \prod_1^m\ev_i^*(\gamma_i)\prod_1^n\widetilde{\ev}_i^*(\widetilde{\gamma}_i), \]
which is an element in $\QQ[\Alb(X)]$ with support on the torsor $T^\delta(L,\beta)$. As the computations in the elliptic case from Section \ref{sec-computation-local-invariants} illustrate, distinct correlators may provide different Gromov--Witten invariants, so that the refinement is non-trivial.

\subsection{Deformation invariants and relations of correlated classes}
Let $(\X,\L)\xrightarrow{p} B$ be a family of smooth projective varieties over $B$ and let $\mathcal L\in \Pic^0(\mathcal X\slash B)$ (in fact log smoothness is sufficient).
Let $\Y=\mathbb P(\O\oplus\L)\to B$ the associated family of $\mathbb P^1$-bundles, with the natural boundary structure $\D=\mathcal D^- +\mathcal D^-$.
Fix $\beta$ an effective curve class on $\X$ relative to~$B$ and
fix contact order $\bfw=(w_1,\dots,w_n)$ of $(q_1,\dots,q_n)$ with the boundary.
As in the previous section, the results of \cite{abramovich2014stable,gross2013logarithmic,maulik2023logarithmic} we have that the moduli stack
$\M_B:=\M_{g, m}(\mathcal Y|\D,\beta_b, \bfw)$
parametrizing diagrams of logarithmic maps
\[
\begin{tikzcd}
C\ar[r,"F"]\ar[d]\ar[d] & (\mathcal Y|\mathcal{D})\ar[d]\\
S\ar[r] & B
\end{tikzcd}
\]
is a Deligne--Mumford stack, proper over $B$ end endowed with a perfect obstruction theory $R\pi_*F^*T^{\log}_{\mathcal Y\slash B}$ relative to $B$.

As for the case of absolute stable maps,
for $B$ regular and connected, the basic properties of the virtual class construction \cite{behrend1997intrinsic,manolache2012virtual} ensure that
the degree of the class $\ev^*\gamma\cap\vir{\M_b}$ does not depend on $b\in B$; we also refer the reader to \cite[Appendix~A]{mandel2020descendant} and references therein for more details on the deformation invariance.

As $\M_B$ is a special case of Corollary~\ref{prop:logmapsdegeneration}, we still have the description of the moduli space in terms of logarithmic trivializations. As in the case of $B=\operatorname{Spec}(\mathbb C)$, we have the map
\[\kappa_B\colon\ \M_B\to\Alb(\mathcal X\slash B)\]
defined by
\[\bigl(\C\xrightarrow{f}\X,(p_j)_{j=1}^m,(q_i)_{i=1}^n,f^*\L\cong\O_{\C}(\alpha)\bigr)\rightarrow \sum_{i=1}^n w_i a_{\X\slash B}(f\circ q_i)\]
which is constant with valued $\varphi_{\beta_b}(\L_b)$ on fibers of $B$. In particular, if $\delta$ divides $\gcd(w_i)$, we get a refinement morphism
\[\kappa^\delta_B\colon\ \M_B\longrightarrow T^\delta_B(\X,\L)\subset \Alb^0(\X/B),
\]
where $T^\delta_B(\X,\L)$ is a the torsor under the $B$-group $\operatorname{Tor}_{\delta}\bigl(\Alb^0(\mathcal X\slash B)\bigr)$ of $\delta$-roots of $\varphi_{\beta_b}(\L_b)$ over~$B$. This torsor is not necessarily trivial, but we can always choose $B'\to B$ \'etale such that $T^\delta_B(\X,\L)\times_B B'\cong T^\delta_{B'}(\X',\L')$ trivializes.

Thus, working \'etale locally on $B$ we can always assume that
$T^\delta_B(\X,\L)$ has a section over~$B$, i.e., it is trivial. Then $\M_B$ split into connected component $\M_B^\theta$ indexed by the sections $\theta\colon B \to T^\delta_B(\X,\L)$.

\begin{theo}\label{thm:definvariants}
The correlated Gromov--Witten invariants are deformation invariant, i.e., let~${\X\to B}$ and $\M_B\to B$ as before and assume that $T^\delta_B(\X,\L)$ has a trivializing section $\theta$. Then~for any such section the degree of classes \smash{$\ev^*\gamma\cap\vir{\M_b^{\theta(b)}}$} for $\gamma\in\text{H}^*(\Y^{n+m})$ does not depend on~${b\in B}$.
\end{theo}

\begin{proof}
 Since \smash{$\M_B^\theta\to \M_B$} is open, the perfect obstruction theory \smash{$R\pi_*F^*T^{\log}_{\mathcal Y\slash B}$} relative to $B$ restricts to $\M_B^\theta$. The results then follows once again from the propertied of the virtual class construction \cite{behrend1997intrinsic,manolache2012virtual} as explained for example in \cite[Appendix~A]{mandel2020descendant}.
\end{proof}

By choosing suitable families of deformations of $\Y\to B$, we can exploit Theorem~\ref{thm:definvariants} to find non-trivial identities among the correlated classes.

The morphism $\varphi_\beta\colon\Pic^0(X)\to\Alb(X)$ induces a morphism between their $\delta$-torsion elements. We thus have the subgroup $\varphi_\beta\bigl(\Tor_\delta\bigl(\Pic^0(X)\bigr)\bigr)\subset\Tor_\delta(\Alb(X))$. In particular, the support $T^\delta(L,\beta)$ of the correlated Gromov--Witten invariants is stable by the action of $\varphi_\beta\bigl(\Tor_\delta\bigl(\Pic^0(X)\bigr)\bigr)$. Furthermore, we have the following.

\begin{theo}\label{theo-torsor-invariance}
 The correlated invariants \smash{$\langle\langle\gamma_1,\dots,\gamma_m,\widetilde{\gamma}_1,\dots,\widetilde{\gamma}_n\rangle\rangle^\delta_{g,\beta,\bfw}$} are invariant under the action of $\varphi_\beta\bigl(\Tor_\delta\bigl(\Pic^0(X)\bigr)\bigr)$.
\end{theo}

\begin{proof}
We consider the family over $B=\Pic^0(X)$ induced by the universal line bundle: the fiber over $L\in\Pic^0(X)$ is $\PP(\O\oplus L)$. Thus, in this case, the torsor appearing before Theorem~\ref{thm:definvariants}~is
\[ T^\delta_B(\X,\L) = \big\{ (L,\theta)\in\Pic^0(X)\times\Alb(X) \text{ s.t. } \delta\theta=\varphi_\beta(L)\in\Alb(X)\big\}.
 \]
To get a section, we use the base change $\mu_\delta\colon L\mapsto L^{\otimes\delta}$ from $\Pic^0(X)$ to itself. We thus have
\[ \mu_\delta^*T^\delta_B(\X,\L) = \big\{ (L,\theta)\in\Pic^0(X)\times\Alb(X) \text{ s.t. } \delta\theta=\varphi_\beta\bigl(L^{\otimes\delta}\bigr)\in\Alb(X)\big\}. \]
This torsor now has a section provided by $L\mapsto\theta_0(L)=\varphi_\beta(L)$.
Deformation invariance tells us that for any choice of $R$, the correlators $\theta_0(L)$ and $\theta_0(L\otimes R)$ provide the same invariants. However, the latter are correlators for the a priori distinct $\PP^1$-bundles associated to $L^{\otimes\delta}$ and~${L^{\otimes\delta}\otimes R^{\otimes\delta}}$. Assuming $R$ is of $\delta$-torsion, these $\PP^1$-bundles are the same. Furthermore, the correlators differ by
\[ \theta_0(L\otimes R)-\theta_0(L) = \varphi_\beta(L\otimes R)-\varphi(L)=\varphi_\beta(R), \]
yielding the result.
\end{proof}

Finally, we have the following relations between the classes corresponding to different levels of refinement.

\begin{lem}\label{lem-unrefinement}
The multiplication by $d$ in $\Alb(X)$ induces a morphism of $\QQ[\Alb(X)]$ denoted by~$\m{d}$. The different full refined virtual classes are related as follows:
\[ \m{\frac{\delta}{\delta'}}\bigl(\fvir{\M}{\delta}\bigr) = \fvir{\M}{\delta'}. \]
\end{lem}

\begin{proof}
 The relation merely comes from the fact that multiplication by $\delta/\delta'$ provides a surjection from $\delta$-roots of $\varphi_\beta(L)$ to $\delta'$-roots of $\varphi_\beta(L)$, so that if $\theta'\in T^{\delta'}(L,\beta)$, then we have
 \[ \M^{\theta'}=\bigsqcup_{(\delta/\delta')\theta=\theta'}\M^{\theta}. \tag*{\qed}
 \]\renewcommand{\qed}{}
\end{proof}

\section{Refined decomposition formula}
Our goal in this section is to give a refinement of the degeneration formula \cite{chen2013degeneration,li2002degeneration} keeping track of our refinement.
 \subsection{Recollection on the degeneration formula}

We consider $\Y_t\xrightarrow{p}\mathbb A^1$ a semi-stable degeneration of a $\mathbb P^1$-bundle $\mathbb P_{\X_{t\neq 0}}(\L\oplus\O)\to\mathbb A^1_{t\neq 0}$ coming from a subdivided tropical line as in the setting of Corollary~\ref{prop:logmapsdegeneration}.

As observed above, that is for example the case when $\Y_t\xrightarrow{p}\mathbb A^1$ is constructed starting from a family of $\mathbb P^1$-bundles over $\X\to\mathbb A^1$ by successively blowing up the boundary divisors on the central fiber.

We recall some notation; write the $N$ components of the central fiber as follows:
\[\Y_0=\bigl(Y_1, D_1^\pm\bigr)\cup \bigl(Y_2, D_2^\pm\bigr)\cup\cdots\cup \bigl(Y_N, D_N^\pm\bigr),\]
where $Y_i$ and $Y_{i+1}$ are gluing along $D_i^+\cong D_{i+1}^-$.

Endowing $\mathbb A^1$ with its toric log structure and $\Y_t$ with the divisorial log structure coming from~${\Y_0+\mathcal D_1^-+\mathcal D_N^+}$, $p$ is logarithmically smooth. Let $C(\mathsf{Y})\xrightarrow{\mathsf{p}} \mathbb R_{\geq 0}$ denote the tropicalization. Then the fiber over any point of $\mathbb R_{>0}$, which we denote by $\mathsf{Y}$, is a subdivided real line with $N$ \emph{vertices} corresponding to the irreducible components of $\Y_0$ and two \emph{legs} corresponding to the divisors~$D_1^-$, $D_N^+$.

Since $\Y_t\xrightarrow{p}\mathbb A^1$ is logarithmically smooth, the moduli space $\M(\Y_t\slash\mathbb A^1)$ of (expanded) logarithmic stable maps to $\Y_t\slash\mathbb A^1$ admits a perfect obstruction theory relative to $\mathbb A^1$ \cite{chen2013degeneration,gross2013logarithmic,li2001stable,maulik2023logarithmic}. In particular, there is a well defined virtual class \smash{$\vir{\M_0}$} for the moduli space of logarithmic stable maps to the central fiber $\Y_0$ and moreover, as recalled above, the logarithmic invariants defined integrating
\smash{$\vir{\M_0}$} and \smash{$\vir{\M_t}$} coincide. This allows one to compute the invariants for maps to the smooth fiber on the degenerate fiber.

The degeneration (or decomposition) formula \cite{chen2013degeneration,kim2018degeneration,li2002degeneration} expresses the virtual fundamental class of the central fiber \smash{$\vir{\M_0}$} as a sum of virtual classes over a collection of decorated graphs, which we call \textit{degeneration graphs}; these encode combinatorial types of curves in the central fiber $\Y_0$.\footnote{ We warn the reader that this terminology is not standard; in \cite{li2002degeneration} degeneration graphs are called \emph{admissible triples}, in \cite{kim2018degeneration} they are referred to as bipartite decorated graph}

\begin{defi}\label{defi-degeneration graph}
A \emph{degeneration graph} for $\Y_0=\bigcup Y_i$ is a graph $\Gamma$ with the following decorations:
\begin{enumerate}\itemsep=0pt
 \item[(1)] We have a graph map $\phi\colon \Gamma\to\mathsf{Y}$, inducing an orientation on the edges of $\Gamma$;

 \item[(2)] \emph{Vertices} are decorated with: a genus $g_V\in\mathbb Z_{>0}$ and a curve class $\beta_V\in H_2(Y_{\phi(V)},\ZZ)$,

 \item[(3)] Every oriented \emph{edge} $e$ and every \emph{leg} $l$ (also called \emph{unbounded edges}) comes decorated with a weight $w_e, w_l\in\mathbb Z$ such that
 the map $\phi$ is \emph{balanced}.

 \item[(4)] The \emph{flow through a vertex} $V$, defined as the sum of the outgoing positive weights, is equal to $\beta_V\cdot \big[D_{\phi(V)}^\pm\big]$.
\end{enumerate}
In what follows, we say \emph{edge} and use the notation $e$ for both edges and legs, unless the distinction is important.

An automorphism of a degeneration graph $\Gamma$ is an automorphism of the underlying graph compatible with the decorations. We denote by $\operatorname{Aut}(\Gamma)$ the group of automorphisms.
\end{defi}

\begin{remarkk}
 In the above definition, \textit{balancing} means that for every vertex $V$, the sum of incoming weights matches the sum of outgoing weights.
\end{remarkk}

\begin{expl}
In Figure~\ref{fig:degeneration}, we draw an example of $C(\mathsf{Y})\to\mathbb R_{\geq 0}$ (on the left) and a corresponding example of degeneration graph (on the right). We represent both the markings carrying non-trivial contact order, namely \emph{unbounded edges} and the internal marking (only $C_{v_3}$ has an internal marking here). Furthermore, each edge, bounded or not, is decorated with the absolute value of its weight; the sign is then determined by the natural orientation.

Finally, each vertex $v$ is further decorated with the data of the genus $g_v$ of the source curve and a curve class $\beta_v\in H_2(X,\mathbb Z)$.
We represent with a cross the components with $g_v=0$ and $\beta_v$ a multiple of the fiber in~$Y_{\phi(v)}$.

\begin{figure}[t]
	\centering
\includegraphics{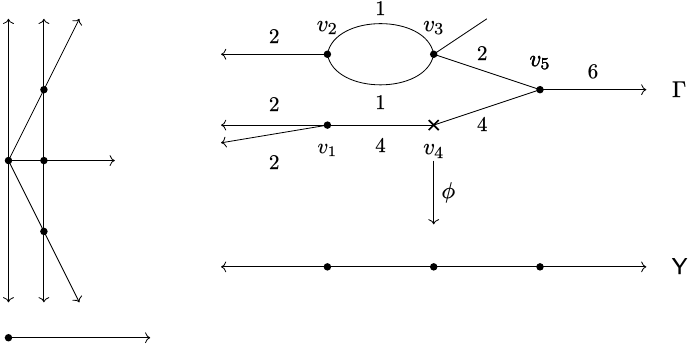}

\caption{Example of degeneration graph.}\label{fig:degeneration}
\end{figure}
\end{expl}

Let $\Gamma$ be a degeneration graph.
Let $V$ be a vertex decorated with genus $g_V$, class $\beta_V$. Let~$n_V$ be the number of adjacent edges in $\Gamma$ with non-zero weight and $m_V$ the number of adjacent legs with trivial weight (called \textit{marked points}). Let $\bfw_e=(w_e)_{e\vdash V}$ be the collection of signed adjacent weights with sign prescribed by the divisors, so that $\sum_{e\vdash V}w_e=0$. We have a moduli space of logarithmic maps associated to each vertex
\smash{$\M_V=\M_{g_V,m_V}(Y_{\phi(V)}|D_{\phi(V)},\beta_V,\bfw_e)$}.

Let $\widetilde{\Gamma}$ be a degeneration graph with a chosen labelling of its bounded edges. For a fixed~$\Gamma$, there are $|E_b(\Gamma)|!\slash|\operatorname{Aut}(\Gamma)|$ such labellings. The cycle version of degeneration formula, as presented by \cite{kim2010logarithmic} (see also \cite{chen2013degeneration,li2002degeneration} for previous versions) states that
\begin{equation}\label{eq:nonrefineddeg}
 \vir{\M_0} = \sum_{\widetilde{\Gamma}} \dfrac{l_{\Gamma}}{|E_b(\Gamma)|!} F_*\phi_{\Gamma}^*\Delta^!\biggl(\prod_V\vir{\M_V}\biggr).
\end{equation}
We recall the notation and refer to \cite{kim2018degeneration} for the proof. We denoted by $l_{\Gamma}=\mathrm{lcm}(\{ w_e\})$; by $\Delta^!$ the Gysin pullback along the diagonal inclusion
\[\bfD:=\prod_e X\to\bfX=\prod_V X^{n_V};\]
by $\phi_{\Gamma}\colon\M_{\widetilde{\Gamma}}\to \bigodot_V \M_V$ the \'etale map of degree \smash{$\frac{\prod_{e\in E_b(\Gamma)} w_e}{l_{\Gamma}}$} to the fiber product
\[
\begin{tikzcd}
\bigodot_V \M_V\ar[r]\ar[d] &\prod_V \M_V\ar[d]\\
\bfD\ar[r,"\Delta"] & \boldsymbol{X}
\end{tikzcd}
\]
corresponding to lifts of the maps in $\bigodot_V \M_V$ to a log stable map to $\Y_0$; by $F\colon \M_{\widetilde{\Gamma}}\to\M_0$ the natural clutching map, shown \cite[Lemma~9.1]{kim2010logarithmic} to have virtual degree \smash{$\frac{|E_b(\Gamma)|!}{l_{\Gamma}}$}.

As also explained in \cite{kim2018degeneration}, and in full details in \cite[Chapter~19]{fulton2013intersection}, via the cycle map \smash{$\operatorname{A}_*\xrightarrow{\text{cl}}\operatorname{H}_{\bullet}$}, we can reduce to perform the computations of the invariance in (Borel--Moore) homology rather then Chow homology. In particular, by \cite[Section~19.1, Theorem~19.2]{fulton2013intersection}
\[\text{cl}\biggl(\Delta^!\biggl(\prod_V\vir{\M_V}\biggr)\biggr)=
\text{cl}\biggl(\biggl(\prod_V\vir{\M_V}\biggr)\biggr)\cap \ev^*\bfD^{\vee}\]
in $H_{\bullet}(\prod_V\vir{\M_V})$
where $\bfD^{\vee}\in H^{\bullet}(\bfX)$ is the Poincar\'e dual of the class of the diagonal.
Then using push-pull formula along the following diagram:
\[
\begin{tikzcd}
\M_{\widetilde{\Gamma}}\ar[r,"F"]\ar[d,"\phi_{\Gamma}"] & \M_0\ar[dr,"ev"]\\
\bigodot_V \M_V\ar[r] & \prod_V\M_V\ar[r,"ev"] & \Y_0
\end{tikzcd}
\]
and the identity above, we then obtain as an immediate corollary from \eqref{eq:nonrefineddeg} the following numerical version of the degeneration formula.
Using the K\"unneth decomposition of $D^\vee$, we can split the cap product over the vertices and get
\begin{equation}
\label{eq-numerical-decomposition-formula}
\left\langle \prod_{i=1}^{n}{\gamma_{i}}\right\rangle_{g,\beta}^{\M_{0}} = \sum_{\Gamma}\frac{\prod_{e}w_e}{|\operatorname{Aut}(\Gamma)|}\int_{\prod\vir{\M_V}} \prod \ev_i^*\gamma_{i}\cup \ev^* D^\vee.
\end{equation}

\subsection{Toward the refined degeneration formula}

For $\Y_t\to \mathbb A^1$, a degeneration coming from a divided tropical line, the interpretation of the~moduli space given in Corollary~\ref{prop:logmapsdegeneration} allows us to define the refinement morphism
$\kappa_{\AA^1} \colon \M\bigl(\Y_t\slash\mathbb A^1\bigr)\to \operatorname{Alb}\bigl(\X\slash\mathbb A^1\bigr)$
as the Albanese evaluation, i.e.,
\[\bigl(\C\xrightarrow{f}\X,(p_j)_{j=1}^m,(q_i)_{i=1}^n,f^*\L\cong\O_{\C}(\alpha)\bigr)\rightarrow \sum_{i=1}^n w_i a_{\X\slash \mathbb A^1}(f\circ q_i).\]

The same argument given in the second part of Proposition~\ref{prop-image-of-L-is-constant} to compute the $\kappa(s)$ in the case of a semi-stable curve
allows us to show that
$\kappa_0\colon \M(\Y_0\slash \operatorname{Spec}(\mathbb N\to\mathbb C))\to \operatorname{Alb}(\X_0)$
is the constant morphism $\varphi_{\beta_0}(\L_0)$. Therefore, as in the case of smooth one parameter families we have that if $\delta$ divides $\gcd(w_i)$ there is a morphism
\smash{$
\kappa_{\AA^1}^\delta\colon \M\bigl(\Y_t\slash\mathbb A^1\bigr)\to T^\delta_{\mathbb A^1}(\X,\L)$}.
Up to passing to an \'etale neighbourhood of $0\in\mathbb A^1$, we can assume that $T^\delta_{\mathbb A^1}(\X,\L)$ trivializes.
In particular, we obtain a decomposition of $\M_0=\M(\Y_0\slash \operatorname{Spec}(\mathbb N\to\mathbb C))$ into connected components indexed by~${\theta\in T^\delta_{\mathbb A^1}(\X,\L)\rvert_0=T^\delta(\X_0,\L_0)}$
\[
\vir{\M_0}=\sum_{\theta\in T^\delta(\X_0,\L_0)} \vir{\M_0^{\theta}}.
\]

Our goal is now to prove a decomposition formula for the components $\vir{\M_0^{\theta}}$.

Let $\Gamma$ be a degeneration graph; for each vertex $V$ let $\delta_V$ is the g.c.d.\ of the weights of the edges adjacent to $V$. Then, as discussed in Section~\ref{sec:correlationrefinement}, we have a refinement
\[ \vir{\M_V}=\sum_{\theta_V\in T^{\delta_V}}\vir{\M_V^{\theta_V}},\]
where $T^{\delta_V}$ is the torsor providing the refinement for $Y_{\phi(V)}$. We aim to describe the components~\smash{$\vir{\M_0^{\theta}}$} in terms of the refined components~\smash{$\M_V^{\theta_V}$}.

\begin{remarkk}
Clearly several choices of $\btheta=(\theta_V)$ may contribute to a common $\M_t^{\theta}$ on the smooth fiber. On the other hand, as we will see shortly, a fixed collection $\btheta=(\theta_V)$ may contribute to several distinct $\M_0^\theta$ and should thus not be considered as a further refinement of the $\M_0^\theta$.
\end{remarkk}

For this reason, the best way to express the refined degeneration formula is in term of \textit{full correlated virtual class} (see Definition~\ref{def:fullrefined}):
\[ \fvir{\M}{\delta} = \sum_{\theta\in T^\delta(L,\beta)} \vir{\M^\theta}\cdot (\theta) \in H_\bullet(\M,\QQ)\otimes\QQ[\Alb(X)].\]

In the group algebra $\QQ[\Alb(X)]$, we also have the division operators defined as follows on generators
\[ \d{\frac{1}{d}}(\theta) = \frac{1}{d^{2r}}\sum_{d\theta'\equiv\theta}(\theta'), \]
where $r$ is the rank of $\Alb(X)$, which is thus a real torus of dimension $2r$. In other words, $\d{\frac{1}{d}}$~maps an element to the average of its $d$-roots.

 \subsection{General degeneration formula}

Let $\Gamma$ be a degeneration graph and $V$ one of its vertices. Let $\delta_\Gamma$ be the g.c.d. of the diagram, i.e., the gcd of the edge weights. Let $n_V$ be its valency, by which we mean the number of adjacent edges and legs, and $m_V$ the number of marked points with trivial contact order. We have the evaluation map
$\ev_V \colon \M_V \to X^{n_V}$
and the Albanese map
\[ a_V\colon\ (x_e)_{e\vdash V}\in X^{n_V}\longmapsto \sum_{e\vdash V}\frac{w_e}{\delta_V}a_X(x_e)\in\Alb(X). \]
The refinement morphism $\kappa_V$ is the composition $a_V\circ\ev_V$.

By definition, the evaluation map $\ev_V$ restricted to the connected component $\M_V^{\theta_V}$ takes values in
$\H_V^{\theta_V}=a_V^{-1}(\theta_V)$.

Let $\btheta=(\theta_V)$ be a vector of correlators indexed by vertices; we denote
\[\boldsymbol{\M}^\btheta=\prod_V \M_V^{\theta_V},\qquad \boldsymbol{\H}^\btheta=\prod_V \H_V^{\theta_V}.\]

Furthermore, we denote by $\widetilde{\ev}_V$ the evaluation map restricted to $\M_V^{\theta_V}$ (omitting $\theta_V$ in the notation) and by $\widetilde{\ev}=\prod_V \widetilde{\ev}_V$, still omitting $\btheta$ in the notation.

As before, we have $\boldsymbol{X}=\prod_V X^{n_V}$ and a natural inclusion $\iota\colon \boldsymbol{\H}^\btheta\hookrightarrow\boldsymbol{X}$ the inclusion, so that we have $\ev=\iota\circ\widetilde{\ev}$. All the data sits in the following diagram:
\[ \begin{tikzcd}
\boldsymbol{\M}^\btheta \arrow[d,"\widetilde{\ev}"] \arrow[rd,"\ev"] & \bfD \arrow[d] & \\
\boldsymbol{\H}^\btheta \arrow[r,"\iota"] & \bfX \arrow[d,"a_{\bfw/\delta}"]\arrow[r,"\mathbf{a}"] & \Alb(X)^{|\V(\Gamma)|}. \\
 & \Alb(X) &
\end{tikzcd} \]

As above, $\bfD$ denotes the diagonal inclusion for each pair of coordinates corresponding to a~bounded edge of $\Gamma$. The components of $\mathbf{a}$ are the $a_V$, while $a_{\bfw/\delta}$ is the morphism giving the global refinement, i.e.,
$a_{\bfw/\delta}(x)=\sum_{e \text{leg}} \frac{w_e}{\delta}a_X(x_e)$.

Clearly, the evaluation map $\prod\M_V\xrightarrow{\ev}\bfX$ may be factored via $\bigsqcup \boldsymbol{\H}^\btheta$ and refined by $\btheta$
\[ \bigsqcup_\btheta \boldsymbol{\M}^\btheta \xrightarrow{\bigsqcup_\btheta\widetilde{\ev}}\bigsqcup_\btheta \boldsymbol{\H}^\btheta \hookrightarrow\bfX .\]
Intersecting the diagonal $\bfD$ with the codomain, the latter also gets refined into \smash{$\bigsqcup_\btheta\bfD^\btheta$}, where we set $\bfD^\btheta=\bfD\cap\boldsymbol{\H}^\btheta$. Therefore, the fiber product $\bigodot\M_V$ also gets refined into $\bigsqcup_\btheta (\bigodot\M_V)^{\btheta}$, and we have the following cartesian diagram:
\[
\begin{tikzcd}
\bigodot_V \M_V\cong \bigsqcup_{\btheta} (\bigodot\M_V)^{\btheta}\ar[r]\ar[d] &\prod_V \M_V=\bigsqcup_{\btheta} \boldsymbol{\M}^\btheta . \ar[d,"\widetilde{\ev}"]\\
\bigsqcup_{\btheta}\bfD^\btheta\ar[r]\ar[d] & \bigsqcup_{\btheta} \boldsymbol{\H}^\btheta\ar[d]\\
\bfD\ar[r,"\Delta"] & \boldsymbol{X}
\end{tikzcd}
\]

The key is that each $\bfD^\btheta$ is itself disconnected and the various distinct components can contribute to distinct $\M_0^\theta$. This contribution is controlled by the following lemma.

\begin{lem}\label{lem-values-varphi-Q-span}
 Let $x\in\bfD^\btheta$. Then, it satisfies
 \[ \frac{\delta}{\delta_\Gamma}a_{\bfw/\delta}(x) = \sum_V\frac{\delta_V}{\delta_\Gamma}\theta_V. \]
\end{lem}

\begin{proof}
Let $x=\bigl(x_e^V\bigr)\in \bfD^\btheta$. We have the following:
\begin{align*}
\frac{\delta}{\delta_\Gamma}a_{\bfw/\delta}(x)&=\frac{\delta}{\delta_\Gamma}\sum_{e \text{leg}} \frac{w_e}{\delta}a_X(x_e)=\sum_{e \text{leg}} \frac{w_e}{\delta_\Gamma}a_X(x_e) = \sum_V \sum_{e\vdash V}\frac{w_e}{\delta_\Gamma}a_X(x_{V,e}) \\
&= \sum_V \frac{\delta_V}{\delta_\Gamma} a_V(x_V)
= \sum_V \frac{\delta_V}{\delta_\Gamma}\theta_V,
\end{align*}
where the third equality follows from the balancing condition that has to be satisfied by those maps that glue to a map to $\Y_0$.
\end{proof}

\begin{remarkk}
 The proof relies on the fact that $a_{\bfw/\delta}$ belongs to the $\QQ$-span of the $a_V$: we have the identity \smash{$\frac{\delta}{\delta_\Gamma}a_{\bfw/\delta}=\sum\frac{\delta_V}{\delta_\Gamma}a_V$}. We may assume that $\delta_\Gamma|\delta$ by choosing $\delta$ to be the g.c.d. of the tangency orders. For other $\delta$, we always can unrefine using Lemma \ref{lem-unrefinement}.
\end{remarkk}

In particular, \smash{$a_{\bfw/\delta}\bigl(\bfD^\btheta\bigr)$} can take a discrete set of values, indexed by the torsor of \smash{$(\delta/\delta_\Gamma)$}-roots of \smash{$\sum\frac{\delta_V}{\delta_\Gamma}\theta_V$}; these values index disjoint components of
$\bfD^\btheta$. Denote by $\bfD^{\btheta,\theta}$ be the component of~$\bfD^\btheta$ corresponding to \smash{$\theta \in a_{\bfw/\delta}\bigl(\bfD^\btheta\bigr)$} and let $\Delta^{\btheta,\theta}$ be the restriction of $\Delta$ to this component then
\smash{$\Delta^{\btheta} =\sum_{\theta}\Delta^{\btheta,\theta}$}.
If $L$ is chosen generically, we can assume that $\boldsymbol{\H}^\btheta$ and $\bfD^\btheta$ are manifolds since they are defined transversally.
In particular, $\Delta^{\btheta}$ is a regular embedding.

%$\Delta\rvert{_\Delta^{\btheta}}$ are regular embeddings
\begin{theo}\label{theo-refined-decomposition}
 In the notation of this section, we have a refined decomposition
 \[ \vir{\M_0^\theta} = \sum_{\widetilde{\Gamma}} \sum_{\btheta} \dfrac{l_{\Gamma}}{|E_b(\Gamma)|!} F_*\phi_{\Gamma}^*\Delta^{{\btheta,\theta},!}\biggl(\prod_V \vir{\M_V^{\theta_V}}\biggr), \]
 where the second sum is over the $\btheta=(\theta_V)$ such that
$\sum\frac{\delta_V}{\delta_\Gamma}\theta_V = \frac{\delta}{\delta_\Gamma}\theta$.
\end{theo}

\begin{proof}
The proof follows from the usual decomposition formula and Lemma \ref{lem-values-varphi-Q-span}. The compatibility of the Gysin pullback from \cite[Chapter~6]{fulton2013intersection} ensures that we may replace $\Delta$ by the sum over the refined diagonals $\Delta^{\btheta,\theta}$.

To get the $\theta$-part of the decomposition for a given $\theta$, we only take the diagonals involving the chosen $\theta$. Following Lemma \ref{lem-values-varphi-Q-span}, this requires to sum over the $\btheta$ satisfying $\sum\frac{\delta_V}{\delta_\Gamma}\theta_V=\frac{\delta}{\delta_\Gamma}\theta$, and we obtain the desired identity in $ H_\bullet(\M,\QQ)\otimes\QQ[\Alb(X)]$.
\end{proof}

In order to extract a numerical version of the refined decomposition formula which can actually be used to compute the invariant, we need
$\deg \bigl(\Delta^{\btheta,\theta,!}\alpha\cap \ev^*\gamma\bigr)$
for $\alpha\!\in\! H_\bullet\bigl(\M^{\btheta},\QQ\bigr)$. This~means that we need expressions for the Poincar\'e dual classes \smash{$\bigl(\bfD^{\btheta,\theta}\bigr)^{\vee}$}.

In the unrefined setting, an explicit expression of the class $\bfD^\vee$ Poincar\'e dual to the diagonal is provided by the K\"unneth decomposition, which allows us to write the class $\bfD^\vee$ in terms of a~basis of $H^*(X,\mathbb Q)$. This provides a decomposition of $\iota^*\bfD^\vee$ but not for the $\bigl(\bfD^{\btheta,\theta}\bigr)^{\vee}$. Furthermore, we do not a priori know much of the cohomology of $\boldsymbol{\H}^\btheta$. Therefore, this task may be especially hard if the classes $\bigl(\bfD^{\btheta,\theta}\bigr)^{\vee}$ are not pulled back from $\mathbf{X}$. %However, we have the following lemma.

 \subsection{Degeneration in the elliptic case}
We first prove the following lemma.

\begin{lem}\label{lem-cobordant-components}
 The diagonal components $\bfD^{\btheta,\theta}$ indexed by \smash{$\theta\in a_{\bfw/\delta}\bigl(\bfD^\btheta\bigr)$} are cobordant in $\bfD$ $($and thus in $\bfX)$. In other words, $\iota_*\big[\bfD^{\btheta,\theta}\big]$ does not depend on the choice of $\theta\in a_{\bfw/\delta}\bigl(\bfD^\btheta\bigr)$.
\end{lem}

\begin{proof}
Let $\theta\in a_{\bfw/\delta}\bigl(\bfD^\btheta\bigr)$ and $\boldsymbol{\gamma}(t)\in\Alb(X)^{|\V(\Gamma)|}$ be a loop based at $\btheta$. The idea is to use the loop $\boldsymbol{\gamma}(t)$ to construct a cobordism between the different components $\bfD^{\btheta,\theta}$. We claim that there is a unique path $\gamma(t)\in\Alb(X)$ such that $\gamma(0)=\theta$ and in $\Alb(X)$ we have
\[ \frac{\delta}{\delta_\Gamma}\gamma(t)=\sum\frac{\delta_V}{\delta_\Gamma}\boldsymbol{\gamma}_V(t). \]
To construct it, let $\widetilde{\theta}$ be a lift of $\theta$ in the universal cover $H^0(X,\Omega_X)^*$ of $\Alb(X)$, and consider the loop $\sum_V \frac{\delta_V}{\delta_\Gamma}\boldsymbol{\gamma}_V(t)$ in $\Alb(X)$. We can lift the loop to a path in $H^0(X,\Omega_X)^*$ starting at~$\frac{\delta}{\delta_\Gamma}\widetilde{\theta}$. This path can be written as follows
\smash{$t\mapsto \frac{\delta}{\delta_\Gamma}\widetilde{\theta}+\rho(t)$},
where $\rho(0)=0$ and $\rho(1)\in H_1(X,\ZZ)\subset H^0(X,\Omega_X)^*$. The path $\rho$ is determined by $\boldsymbol{\gamma}$. It corresponds to a loop in $\Alb(X)$, and its class $[\rho]=\rho(1)\in H_1(X,\ZZ)\subset H^0(X,\Omega_X)^*$ satisfies $[\rho]=\Sigma_*(\boldsymbol{\gamma})$, where $\Sigma_*\colon H_1(\Alb(X)^{|\V(\Gamma)|},\ZZ)\to H_1(\Alb(X),\ZZ)$ is induced by \smash{$\Sigma\colon(\theta'_V)\mapsto\sum\frac{\delta_V}{\delta_\Gamma}\theta'_V$}.
Finally, we divide by $\delta/\delta_\Gamma$: the path we care about is the image in $\Alb(X)$ of \smash{$t\mapsto \widetilde{\theta}+\frac{1}{\delta/\delta_\Gamma}\rho(t)$}.
In particular, $\gamma(t)\in\Alb(X)$ may not be a~loop.

Starting with a generic path $\boldsymbol{\gamma}(t)$, the path $\gamma(t)\in\Alb(X)$ is also generic. Therefore, its preimage by \smash{$a_{\bfw/\delta}|_{\bfD}\colon \bfD\to\Alb(X)$} provides a cobordism in $\bfD$ between $\bfD^{\btheta,\gamma(0)}$ and $\bfD^{\btheta,\gamma(1)}$.

To finish the proof, we need to prove that $\gamma(1)$ may take any value in $a_{\bfw/\delta}\bigl(\bfD^\btheta\bigr)$. Equivalently, we need to show that any element of $\Tor_{\delta/\delta_\Gamma}(\Alb(X))$ has a lift in $H^0(X,\Omega_X)^*$ of the form~\smash{$\frac{1}{\delta/\delta_\Gamma}\rho(1)$}, with $\rho$ as above. This amounts to the surjectivity of $\Sigma_*$. The latter is ensured by the fact that the g.c.d. of $\delta_V/\delta_\Gamma$ is $1$, finishing the proof.
\end{proof}

\begin{remarkk}
 The key point in the proof is the surjectivity. Otherwise, the elements in the torsor would split in different classes modulo the image of the morphism.
\end{remarkk}

\begin{remarkk}
 This lemma is easier to prove in the elliptic case since all maps are actually group morphisms, so that the various components are in fact parallel subtori.
\end{remarkk}

Now, let us work under the additional assumption that the $\iota^*\colon H^\bullet(X^{n_V},\QQ)\to H^\bullet\bigl(\H_V^{\theta_V},\QQ\bigr)$ are surjective. This hypothesis is satisfied when $X$ is an elliptic curve, since each \smash{$\H_V^{\theta_V}$} is now a~subtorus of $X^{n_V}$.
Using Lemma \ref{lem-cobordant-components}, we can now prove the following.

\begin{lem}\label{lem-equal-intersection-numbers}
Given $\gamma\in H^\bullet(\bfX,\QQ)$ and $[\N]$ a cycle class in $H_\bullet\bigl(\boldsymbol{\M}^\btheta,\QQ\bigr)$, the intersection numbers
\smash{$\widetilde{\ev}_*[\N]\cap\bigl(\bigl(\bfD^{\btheta,\theta}\bigr)^{\vee}\cup\iota^*\gamma\bigr) $}
do not depend on $\theta\in a_{\bfw/\delta}\bigl(\bfD^\btheta\bigr)$:
\end{lem}

\begin{proof}
If $L$ is chosen generically, every choice of $\btheta$ is also generic and $\H_V^{\theta_V}$ is thus a submanifold of $X^{n_V}$. In particular, we have Poincar\'e duality. Let $\mu$ be the class Poincar\'e dual to $\widetilde{\ev}_*[\N]$ inside $\boldsymbol{\H}^\btheta$, so that we have
\[ \widetilde{\ev}_*[\N]\cap \bigl(\bigl(\bfD^{\btheta,\theta}\bigr)^{\vee}\cup\iota^*\gamma\bigr) = \big[\boldsymbol{\H}^\btheta\big]\cap\bigl(\bigl(\bfD^{\btheta,\theta}\bigr)^{\vee}\cup\mu\cup\iota^*\gamma\bigr) = \big[\bfD^{\btheta,\theta}\big]\cap(\mu\cup\iota^*\gamma), \]
since $\bigl(\bfD^{\btheta,\theta}\bigr)^{\vee}$ is Poincar\'e dual to $\bfD^{\btheta,\theta}$ inside $\boldsymbol{\H}^\btheta$. By surjectivity, we can write $\mu=\iota^*\widetilde{\mu}$. Moreover, as we compute an intersection number, we may as well compute its push-forward inside $H_0(\bfX,\QQ)$ and use push-pull formula, so that the number we care about is
\[ \iota_*\bigl(\big[\bfD^{\btheta,\theta}\big]\cap\iota^*(\widetilde{\mu}\cup\gamma)\bigr) = \iota_*\big[\bfD^{\btheta,\theta}\big]\cap (\widetilde{\mu}\cup\gamma). \]
As by Lemma \ref{lem-cobordant-components} $\iota_*\big[\bfD^{\btheta,\theta}\big]$ does not depend on $\theta$, we get that the intersection numbers do not depend on the particular choice of $\theta\in a_{\bfw/\delta}\bigl(\bfD^\btheta\bigr)$.
\end{proof}

We may now state the degeneration formula under the surjectivity assumption.

\begin{theo}\label{theo-refined-decomposition-surjective-case}
Consider a one parameter family degeneration $\Y_t$ of $\PP^1$-bundles with common base $X$ whose central fiber is a union of $\PP^1$-bundle glued over their boundary divisors. Under the surjectivity assumption, we have the following decomposition of the virtual classes:
\[ \fvir{\M_0}{\delta}
\equiv \sum_\Gamma\sum_\btheta \frac{\prod w_e}{|\operatorname{Aut}(\Gamma)|}\biggl(\Delta^!\prod_V\vir{\M_V^{\theta_V}} \biggr) \d{\frac{1}{\delta/\delta_\Gamma}}\biggl(\sum\frac{\delta_V}{\delta_\Gamma}\theta_V\biggr), \]
where $\Delta\colon\bfD\to\bfX$ is the usual diagonal inclusion.
\end{theo}

The ``$\equiv$" means that the equality is true when we compute the intersection with the pullback of a cohomology class by the evaluation map.

\begin{proof}
 We just rebrand Theorem~\ref{theo-refined-decomposition} using Lemma \ref{lem-equal-intersection-numbers} to show that both classes have the same intersection numbers with classes of the form $\iota^*\gamma$. Notice that since this is a statement about intersection numbers, as in equation \eqref{eq-numerical-decomposition-formula}, we replace the sum over labeled graphs $\widetilde{\Gamma}$ by a sum over degeneration graphs $\Gamma$, each one having precisely $\frac{|E_b(\Gamma)|!}{|\operatorname{Aut}(\Gamma)}$ labellings. Taking into account the degree of $\phi_\Gamma$ and the virtual degree of $F$ gives the coefficient for each $\Gamma$, which is independent of the refinement.

 Let \smash{$\vir{\boldsymbol{\M}^\btheta}=\prod\vir{\M_V^{\theta_V}}$} and $r$ be the dimension of $\Alb(X)$, so that the cardinality of $\delta/\delta_\Gamma$-torsion elements is $(\delta/\delta_\Gamma)^{2r}$. Intersection numbers in the (co)homology of $\boldsymbol{\H}^\btheta$ between a~class in~$H_\bullet\bigl(\boldsymbol{\H}^\btheta,\QQ\bigr)$, some $\bigl(\bfD^{\btheta,\theta}\bigr)^{\vee}$ and some $\iota^*\gamma$ can be computed in the (co)homology of $\bfX$ by push-pull-formula. By Lemma \ref{lem-equal-intersection-numbers}, the results do not depend on the specific choice of $\theta$. Therefore, they are equal to their average
 \begin{align*}
 \iota_*\bigl(\widetilde{\ev}_*\vir{\boldsymbol{\M}^\btheta}\cap \bigl(\bigl(\bfD^{\btheta,\theta}\bigr)^{\vee}\cup\iota^*\gamma\bigr) \bigr) ={} & \iota_*\left( \widetilde{\ev}_*\vir{\boldsymbol{\M}^\btheta}\cap \left(\frac{1}{(\delta/\delta_\Gamma)^{2r}}\sum_{\theta}\bigl(\bfD^{\btheta,\theta}\bigr)^{\vee}\cup\iota^*\gamma\right) \right) \\
 ={} & \frac{1}{(\delta/\delta_\Gamma)^{2r}} \iota_*\bigl(\widetilde{\ev}_*\vir{\boldsymbol{\M}^\btheta}\cap \iota^*\bfD^\vee\cup\iota^*\gamma \bigr) \\
 = {}& \frac{1}{(\delta/\delta_\Gamma)^{2r}} \ev_*\vir{\boldsymbol{\M}^\btheta}\cap\bigl(\bfD^\vee\cup\gamma\bigr).
 \end{align*}
 As the $\frac{1}{(\delta/\delta_\Gamma)^{2r}}$ are precisely the coefficients appearing through the use of $\dfk\big[\frac{1}{\delta/\delta_\Gamma}\big]$ and its support the $\theta$ we care about, we get the result.
\end{proof}

\begin{coro}
Under the surjectivity assumption, the full correlated virtual classes satisfy the following decomposition formula:
\[ \fvir{\M_0}{\delta}
\equiv \sum_\Gamma \frac{\prod w_e}{|\operatorname{Aut}(\Gamma)|}\d{\frac{1}{\delta/\delta_\Gamma}}\Delta^!\prod_V \fvir{\M_V}{\delta_\Gamma}. \]
\end{coro}

\begin{proof}
We rewrite Theorem~\ref{theo-refined-decomposition-surjective-case} to make appear the full refined classes at the vertex level
	\begin{align*}
	\prod_V\fvir{\M_V}{\delta_\Gamma} = & \prod_V \m{\frac{\delta_V}{\delta_\Gamma}}\biggl(\sum_{\theta_V} \vir{\M_V^{\theta_V}}\cdot (\theta_V) \biggr)
	= \prod_V\sum_{\theta_V} \vir{\M_V^{\theta_V}}\cdot \left(\frac{\delta_V}{\delta_\Gamma}\theta_V\right) \\
	= & \sum_{\btheta} \prod_V\vir{\M_V^{\theta_V}}\cdot \left(\sum\frac{\delta_V}{\delta_\Gamma}\theta_V \right),
	\end{align*}
	yielding the result.
\end{proof}

\begin{remarkk}
The formula tells us that even if each vertex is entitled to a refinement at the level $\delta_V$, to recover the refinement for the global degeneration graph at the level $\delta$, we only need to know the refinements at the level $\delta_\Gamma$. Moreover, due to the presence of $\dfk\big[\frac{1}{\delta/\delta_\Gamma}\big]$, the class associated to the graph $\Gamma$ is invariant under $\Tor_{\delta/\delta_\Gamma}(\Alb(X))$. In particular, if $\delta_\Gamma=1$, all the correlated classes coming from $\Gamma$ yield the same invariants.
\end{remarkk}

\section{Computation of local invariants}
\label{sec-computation-local-invariants}

Our goal in this section is to compute the correlated GW invariants in the case where $X=E$ is an elliptic curve, the genus of the source curve is $g=1$, and with a unique interior point constraint. These \textit{local invariants} will be used in combination with Theorem~\ref{theo-refined-decomposition-surjective-case} to obtain an explicit computation algorithm and derive regularity results for general correlated invariants in Section \ref{sec-regularity}. We provide a much shorter and less technical proof than the one provided in the previous version of this paper.

We consider $Y=\PP(\O\oplus L)$ for $L\in\Pic^0(E)$; fix an homology class $\beta=a[E]\in H_2(E,\ZZ)$ and a vector $\bfw=(w_1,\dots,w_n)$ of tangency orders. The moduli space
$\M(a,\bfw)=\M_{1,1}\bigl(Y|D^\pm,a[E],\bfw\bigr)$ is the moduli space of log stable maps as in the previous section.
For $\delta|\mathrm{gcd}(w_i)$, we saw that it decomposes into components $\M^\theta(a,\bfw)$.
The correlators $\theta$ satisfy
$
\delta\theta\equiv \varphi_{a[E]}(L)=a\lambda$,
where~\smash{$\lambda=\varphi_{[E]}(L)$} is the image of the line bundle $L\in\Pic^0(E)\simeq E$ through the isomorphism~$\varphi_{[E]}$. We denote the torsor by $T^\delta(L,a)$. Up to deformation, we can assume that $L=\O$ and $T^\delta(L,a)$ is the torsor of $\delta$-torsion elements $\Tor_\delta(E)$.

We consider the correlated invariants
\[ \langle\langle \pt_0,1_{w_1},\pt_{w_2},\dots,\pt_{w_n}\rangle\rangle^\delta_{1,\beta,\bfw}
= \int_{[\![\M(a,\bfw)]\!]^\delta}\ev_0^*(\pt)\prod_2^n\ev_i^*(\pt), \]
which is an element of the group algebra $\QQ[E]$ with support on $T^\delta(L,a)$. To provide a compact formula, we introduce some elements of $\QQ[E]$: the average of torsion elements $T_d$ is defined by
\[T_d:=\frac{1}{d^2}\sum_{d\theta\equiv 0}(\theta) \in\QQ[E],\]
which we use to define a correlated version of the sum of divisors function: for $\delta\geqslant 1$, we set
\[ \boldsymbol{\sigma}_\delta(a) = \sum_{k|a}\frac{a}{k}T_{\delta/\gcd(k,\delta)} .\]

\begin{lem}
The function $a\mapsto\boldsymbol{\sigma}_\delta(a)$ is a quasi-modular form for the congruence subgroup $\Gamma_0(\delta)$ with values in $\QQ[E]$.
\end{lem}

\begin{proof}
Since $\boldsymbol{\sigma}_\delta$ has values in the finite dimensional subalgebra generated by the $(T_d)_{d|\delta}$, it suffices to check that its coordinates are quasi-modular. Grouping the $T_d$ together,
\[\boldsymbol{\sigma}_\delta(a) = \sum_{d|\delta}\left(\sum_{\substack{k|a \\ \gcd(k,\delta)=d}}\frac{a}{k}\right)T_{\delta/d}.\]
Writing $k=dk'$ we have for every $d|\delta$
\[\sum_{\substack{k|a \\ d|\gcd(k,\delta)|\delta}}\frac{a}{k} = \sum_{d|k|a}\frac{a}{k} = \sum_{k'|\frac{a}{d}}\frac{a/d}{k} =\sigma(\frac{a}{d}),\]
understood to be $0$ if $d$ does not divide $a$. Performing a M\"obius transformation using the M\"obius function $\mu$, we deduce that
\[\sum_{\substack{k|a \\ \gcd(k,\delta)=d}}\frac{a}{k} = \sum_{d|l|\delta}\mu(\frac{l}{d})\sigma(\frac{a}{l}).\]

The generating series of $\sigma$ is the second Eisenstein series $E_2(\sfq)$, and it is a quasi-modular form for $SL_2(\ZZ)$. Hence, for $d|\delta$, we have
\[\sum_{a=1}^\infty \sigma(\frac{a}{d})\sfq^a = \sum_{a'=1}^\infty \sigma(a')\sfq^{da'} = E_2(\sfq^d),\]
which is thus a quasi-modular form for the congruence subgroup $\Gamma_0(d)\supset\Gamma_0(\delta)$. In particular, so is $\boldsymbol{\sigma}_\delta$ since its coordinates are.
\end{proof}

\begin{theo}\label{theo-expression-local-correlated-invariant}
The full local correlated invariant has the following expression:
\begin{equation}\label{eq-local-invariant}
\langle\langle \pt_0,1_{w_1},\pt_{w_2},\dots,\pt_{w_n}\rangle\rangle^\delta_{1,a[E],\bfw} = a^{n-1}w_1^2 \boldsymbol{\sigma}_\delta(a)
\cdot (\theta_0),
\end{equation}
where $\theta_0$ is any special correlator, i.e., satisfies \smash{$\frac{\delta}{\gcd(a,\delta)}\theta_0 = \frac{a}{\gcd(a,\delta)}\lambda$}.
\end{theo}

Without assuming Theorem~\ref{theo-expression-local-correlated-invariant}, it can be checked from the definition of $\boldsymbol{\sigma}$ that
\begin{equation}\label{eq-sigma-invariance}
 \boldsymbol{\sigma}_\delta(a) = \boldsymbol{\sigma}_\delta(a)T_{\delta/\gcd(a,\delta)},
\end{equation}
so that the right-hand side of \eqref{eq-local-invariant} does not depend on which $\theta_0$ we choose. Assuming instead Theorem~\ref{theo-expression-local-correlated-invariant}, equation \eqref{eq-sigma-invariance} is in fact a consequence from Theorem~\ref{theo-torsor-invariance}, which tells us that the correlated invariant is invariant when multiplying by an element of \[\varphi_{a[E]}(\Tor_\delta(E))=a\cdot\Tor_\delta(E)=\Tor_{a/\gcd(a,\delta)}(E).\]
In particular, it is also invariant by multiplication by $T_{\delta/\gcd(a,\delta)}$.

With the above formulation, the correlated invariant appears as a refined version of the uncorrelated invariant $a^{n-1}w_1^2\sigma(a)$, and the refinement takes the form of $\boldsymbol{\sigma}_\delta$ replacing $\sigma$, the sum of divisors function.

Assume that $L=\O$ and choose $\theta_0=0$, so that the torsor of correlators is actually $\Tor_\delta(E)$. For $\theta\in\Tor_\delta(E)$, its order $\omega(\theta)$ is the smallest positive integer $n$ such that $n\theta\equiv 0$.

\begin{lem}\label{lem-monodromy-invariance}
 Assuming $L=\O$, the correlated invariant \smash{$\langle \pt_0,1_{w_1},\pt_{w_2},\dots,\pt_{w_n}\rangle_{1,a[E],\bfw}^\theta$} only depends on $\theta$ through its order $\omega(\theta)$.
\end{lem}

\begin{proof}
 Write $E$ as some $E_\tau=\CC/\langle1;\tau\rangle$ for some complex number $\tau$ in the Poincar\'e half-plane. For any \smash{$\gamma=\bigl(\begin{smallmatrix}
 a & b \\ c & d \\
 \end{smallmatrix}\bigr)\in {\rm SL}_2(\ZZ)$}, we have an isomorphism between $E_\tau$ and $E_{\gamma\cdot\tau}$ where \smash{$\gamma\cdot\tau=\frac{a\tau+b}{c\tau+d}$}. The isomorphism is actually given by
 \[ \varphi\colon\ (z \modulo 1,\tau) \longmapsto \frac{z}{c\tau+d} \quad \modulo 1,\gamma\cdot\tau. \]
 Choose a path $\tau(t)$ in the Poincar\'e half-plane going from $\tau(0)=\tau$ to $\tau(1)=\gamma\cdot\tau$. The $\delta$-torsion elements in $E_{\tau(t)}$ are deformed continuously along with $\tau(t)$ as they can be written as $u\tau(t)+v$, where $u,v\in\Tor_\delta(\RR/\ZZ)$. The identification $\varphi\colon E_\tau\to E_{\gamma\cdot\tau}$ then induces some monodromy among torsion elements. More precisely, if $u,v\in \Tor_\delta(\RR/\ZZ)$, we have
 \[ \varphi^{-1}\left(u\frac{a\tau+b}{c\tau+d}+v \right)=u(a\tau+b)+v(c\tau+d) = (au+cv)\tau+(bu+dv). \]
 Therefore, we deduce that the values of the correlated invariants for the torsion elements $\theta=u\tau+v$ and $\gamma\cdot\theta=(au+cv)\tau+(bu+dv)$ are the same. As the action of ${\rm SL}_2(\ZZ)$ on the set of torsion elements of the same order in $\Tor_\delta\bigl((\RR/\ZZ)^2\bigr)$ is transitive, we conclude.
\end{proof}

\begin{proof}[Proof of Theorem \ref{theo-expression-local-correlated-invariant}]
Writing $E=\CC/\gen{1,\tau}$, we get a group isomorphism $E\simeq \RR/\ZZ\times\RR/\ZZ$. The projection $\psi\colon E\to \RR/\ZZ$ along the second coordinate (the $\tau$-coordinate) provides a morphism of group algebra $\psi_*\colon\QQ[E]\to\QQ[\RR/\ZZ]$. Denoting by $T'_d=\psi_*(T_d)=\frac{1}{d}\sum_{d\theta'}(\theta')$ the average of $\delta$-torsion elements in $\QQ[\RR/\ZZ]$, $\psi$ induces an isomorphism from the subalgebra spanned by $(T_d)_{d|\delta}$ to the algebra spanned by $(T'_d)_{d|\delta}$.

We assume $L=\O$, so that $\theta_0=0$ is a special correlator. The monodromy invariance from Lemma \ref{lem-monodromy-invariance} ensures that the full correlated invariant belongs to the span of $(T_d)_{d|\delta}$. Therefore, to compute it, it suffices to compute its image by $\psi_*$.

A genus $1$ stable map $C\to E\times\PP^1$ consists in a covering $f\colon C\to E$ of degree $a$ together with a meromorphic function $C\to\PP^1$ with prescribed poles and zeroes. Using the $0$th-marked point as neutral element of $C$ and its image as neutral element of $E$, the set of stable maps meeting the point constraints $(p_0,x_2,\dots,x_n)\in (E\times\PP^1)\times E^{n-1}$ is in bijection with the following set:
\[ \left\{ (f\colon C\to E,y_1,\dots,y_n) \text{ s.t. }\forall\ 2\leqslant i\leqslant n\; f(y_i)=x_i, \ \sum_1^n w_iy_i=0\in C\right\},\]
where $f\colon C\to E$ is a degree $a$ covering and $y_1,\dots,y_n\in C$ are points.

Degree $a$ coverings are in bijection with index $a$ sublattices of $\ZZ+\tau\ZZ\simeq\pi_1(E)$. Each sublattice has a unique basis of the form $e_1=\frac{a}{k}$ and $e_\tau=j+k\tau$, where $k|a$ and $0\leqslant j<\frac{a}{k}$ are integers. In particular, there are $\sigma(a)=\sum_{k|a}\frac{a}{k}$ such coverings. With this notation, $C=\CC/\gen{e_1,e_\tau}$. The $\delta$-torsion elements of $C$ may be written $ue_1+ve_\tau$ where $u,v\in(\frac{1}{\delta}\ZZ)/\ZZ$. Their image in $E$ is still $ue_1+ve_\tau$ but taken this time modulo $1,\tau$.

We compute the contribution of a fixed covering $f\colon C\to E$ to the correlated invariant. We fix such a covering, i.e. a choice of $(k,j)$ with $k|a$ and $0\leqslant j<\frac{a}{k}$. There are $a$ possible choices of $y_i\in C$ for each $x_i\in E$ and $2\leqslant i\leqslant n$. For each choice of $(y_2,\dots,y_n)$, the last point $y_1$ is fixed by the condition $\sum_1^n w_iy_i=0$. Choosing a fixed $y_1$ such that $\sum_1^n\frac{w_i}{\delta}y_i=0$, the other possible $y_1$ are obtained adding some $w_1$-torsion element $z$, thus inducing a bijection between the solutions with given $(f\colon C\to E,y_2,\dots,y_n)$ and the $w_1$-torsion elements. The correlator satisfies $\theta=\sum_1^n \frac{w_i}{\delta}f(y_i)=\frac{w_1}{\delta}f(z)$. Therefore, the count of solutions associated to the covering $f\colon C\to E$ accounting for the correlator is
\[\operatorname{Mult}(f) = a^{n-1}\sum_{z\in \Tor_{w_1}(C)}(\frac{w_1}{\delta}f(z)).\]
Since multiplication by $\frac{w_1}{\delta}$ induces a surjection from $w_1$-torsion to $\delta$-torsion, and each $\delta$-torsion element has precisely $(\frac{w_1}{\delta})^2$ preimages, we get
\[\operatorname{Mult}(f) = a^{n-1}(\frac{w_1}{\delta})^2\sum_{t\in \Tor_\delta(C)}(f(t)).\]
Notice that this element has no reason to belong to the span of $(T_d)_{d|\delta}$. We now apply $\psi_*$, projecting the group generators onto their second coordinate. Since $e_1$ has $\tau$-coordinate equal to $0$, the second coordinate of $f(t)$ takes value in the subgroup generated by the second coordinate of $e_\tau$, which is $\frac{k}{\delta}$ in $(\frac{1}{\delta}\ZZ)/\ZZ$. By Bezout theorem, this subgoup is also generated by $\frac{\gcd(k,\delta)}{\delta}$, and these are exactly the $\delta/\gcd(k,\delta)$-torsion elements. Therefore, we deduce that
\[\psi_*(\operatorname{Mult}(f)) = a^{n-1}w_1^2 T'_{\delta/\gcd(k,\delta)}.\]
In particular, it only depends on $k$. Summing over the possible values of $0\leqslant j<\frac{a}{k}$ and $k$, we get
\[\sum_{f\colon C\to E}\psi_*(\operatorname{Mult}(f)) = a^{n-1}w_1^2\sum_{k|a}\frac{a}{k}T'_{\delta/\gcd(k,\delta)}.\]
Using that $\psi$ induces an isomorphism from the algebra spanned by $(T_d)_{d|\delta}$, we deduce that
\[\sum_{f\colon C\to E}\operatorname{Mult}(f) = a^{n-1}w_1^2\sum_{k|a}\frac{a}{k}T_{\delta/\gcd(k,\delta)}.\]
Stress that it does not imply that each $\operatorname{Mult}(f)$ is a multiple of some $T_d$, and this is actually false.
\end{proof}

\section{Floor diagrams and regularity}
\label{sec-regularity}

Our goal in this section is to obtain regularity results in the flavor of \cite{blomme2022floor} in the refined setting provided by the correlation. To do so, we use the correlated degeneration formula from Theorem~\ref{theo-refined-decomposition-surjective-case} along with the local computation from Theorem~\ref{theo-expression-local-correlated-invariant}. The use of the decomposition formula for a family of $E\times\PP^1$ with central fiber a chain of $E\times\PP^1$ glued along their boundary divisors leads to combinatorial objects called \textit{floor diagrams}. The latter appear in \cite{blomme2022floor} although they are obtained through tropical techniques. Using these floor diagrams, we are able to prove the quasi-modularity in the base direction, and the piecewise polynomiality in the tangency orders. We explain how to adapt the floor diagram to the correlated setting, using coefficients in the group algebra $\QQ[E]$.

\subsection{Floor diagrams and their multiplicities}

We first recall the notion of \textit{floor diagram} in the setting of curves in $E\times\PP^1$, already developed in~\cite{blomme2022floor}.\looseness=-1

\begin{defi}
 A floor diagram $\Dfk$ is the data of a weighted oriented graph with the following properties:
 \begin{enumerate}\itemsep=-0.5pt
 \item[(1)] The graph has three kind of vertices:
 \begin{itemize}\itemsep=-0.5pt
 \item sinks and sources which are univalent vertices, referred as \textit{infinite vertices},
 \item \textit{flat vertices} which are bivalent with one ingoing and one outgoing edge,
 \item \textit{floors} which carry a label $a_V\in\NN$.
 \end{itemize}
 \item[(2)] The set of flat vertices and floors carry a total order compatible with the orientation.
 \item[(3)] The edges have a positive weight $w_e$, such that the weighting makes floors and flat vertices balanced.
 \item[(4)] The complement of all flat vertices is without cycle and each connected component of the complement contains a unique sink or source.
 \end{enumerate}
\end{defi}

For floor diagram, as for curves, we have a notion of genus and class, recalled in the next definition.

\begin{defi}\samepage
 Let $\Dfk$ be a floor diagram.
 \begin{itemize}\itemsep=-0.5pt
 \item The genus $g(\Dfk)$ is its genus as a graph where floors are considered to have genus $1$: $g(\Dfk)=b_1(\Dfk)+|V(\Dfk)|$, where $V(\Dfk)$ is the set of floors.
 \item The class of $\Dfk$ is $a[E]+b\big[\PP^1\big]$ where $a=\sum_V a_V$ is the sum of floors weights, and $b$ is the sum of flows at the sources (or sinks, since the flow through the diagram is preserved by the balancing condition).
 \item The tangency profile $\bfw$ is the collection of $\pm w_e$ for edges adjacent to infinite vertices (with a $+$ for sinks and $-$ for sources).
 \item The $\delta$-gcd of a diagram is the gcd of its edge weights and $\delta$.
 \end{itemize}
\end{defi}

\begin{figure}[th]
 \centering
\includegraphics{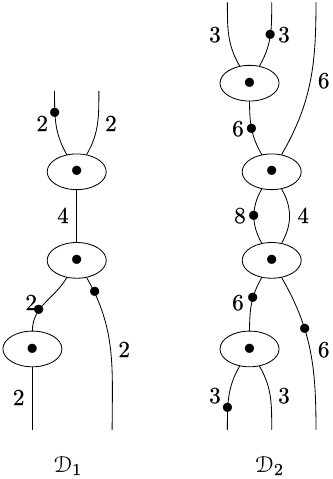}
 \caption{Floor diagrams.}
 \label{fig-expl-diagrams}
\end{figure}

We depict floor diagrams as on Figure \ref{fig-expl-diagrams}, with edges oriented from bottom to top: sources are at the bottom and sinks at the top. Floors are big-shaped with a bullet while flat vertices are bullets on edges.

\begin{expl}
\label{expl-floor-diagram}
 We have two different floor diagrams on Figure \ref{fig-expl-diagrams}. The first diagram has genus $3$, gcd $2$ and tangency profile $(2,2,-2,-2)$. The second diagram has genus $4$ and gcd $1$, although the gcd of the weights at infinity is $3$. The tangency profile is $(3,3,6,-3,-3,-6)$. For each diagram, we can check that all vertices are balanced and the complement of all flat vertices is a~forest.
\end{expl}

Floor diagrams encode some combinatorial types of curves in an expanded degeneration of~${E\times\PP^1}$, which is the gluing of copies of $E\times\PP^1$ along their boundary divisors. As we prove in Theorem~\ref{theo-correspondence-diagrams}, the latter are exactly the graphs appearing in the decomposition formula, up to a~redecoration.
 \begin{itemize}\itemsep=0pt
 \item Sources and sinks correspond to marked points mapped to the boundary divisors of the expanded degeneration.
 \item Flat vertices correspond to genus $0$ components with an internal marked point mapped to a fiber of $E\times\PP^1$.
 \item Floors correspond to genus $1$ components in the class $a_V[E]+b_V\big[\PP^1\big]$ where $b_V$ is the flow through $V$. Its intersection profile with the boundary divisor is prescribed by the adjacent edges.
 \item Edges correspond to nodes in the domain mapped to the singular locus of the expanded degeneration.
 \item Edges relating vertices which are not consecutive in the order on floors and flat vertices should be subdivided. The new bivalent vertices would correspond to components mapped to a fiber but without marked point.
 \end{itemize}

Let $\bfw=(w_1,\dots,w_n)$ be a tangency profile and let $\delta$ be a common divisor. Let $\Dfk$ be a floor diagram and $\delta_\Dfk$ its $\delta$-gcd (the gcd between its edge weights and $\delta$). Recall that we have the division operator $\dfk\big[\frac{1}{\delta/\delta_\Dfk}\big]$ over the group algebra.

\begin{defi}\label{defi-multiplicity-diagram}
 We define the correlated multiplicity of $\Dfk$ to be
 \[ \bfm_\delta(\Dfk) = \d{\frac{1}{\delta/\delta_\Dfk}}\biggl(\prod_V a_V^{n_V-1}\boldsymbol{\sigma}_{\delta_\Dfk}(a_V) \biggr) \prod_{e\in E_b(\Dfk)} w_e\prod_{e\in E_\circ(\Dfk)}w_e^2 \cdot (a\lambda_0) \in \QQ[E], \]
 where the first product is over the floors, the second product over the set of bounded edges~$E_b(\Dfk)$, the last product over the set of possibly unbounded edges not adjacent to a flat vertex $E_\circ(\Dfk)$, and $\delta\lambda_0=\lambda$.
\end{defi}

In the literature (for instance, \cite{blomme2022floor,brugalle2008floor}), the flat vertices are not vertices but they are considered as markings on the graph. Deleting flat vertices and merging the adjacent edges sharing the same weight, the exponent of the weight $w_e$ of a bounded edge (resp.\ end) in the multiplicity would be $2$ (resp.\ $1$) if the edge carries a marking, and $3$ (resp.\ $2$) if it is not, matching the result from \cite{blomme2022floor}.

This multiplicity $\bfm_\delta$ is a refinement of the multiplicity presented in \cite{blomme2022floor}, which is actually the multiplicity $\bfm_1$.

\subsection{Correspondence statement}

We now use the degeneration formula to prove that the count of floor diagrams with the ad-hoc multiplicity yields the correlated GW-invariant. We consider a degeneration $\Y_t$ of $\PP(\O\oplus L)$ into~${n+g-1}$ components and one point constraint in each component. In particular, for a~curve $f\colon C\to\Y_0$, among components of $C$ mapped to a common component of $\Y_0$, exactly one may contain a marked point.

\begin{theo}\label{theo-correspondence-diagrams}
 We have the following equality:
 \[ \gen{\gen{\pt^{n+g-1},1_{w_1},\dots,1_{w_n}}}^\delta_{g,a[E],\bfw} =
 \sum_\Dfk \bfm_\delta(\Dfk), \]
 where we sum over floor diagrams $\Dfk$ of genus $g$, in the class $\beta$, having tangency profile $\bfw$.
\end{theo}

The proof is standard application of the decomposition formula which we recall for sake of completeness.

Up to the factor $\frac{\prod w_e}{|\mathrm{Aut}(\Gamma)|}$ and the operator $\dfk\big[\frac{1}{\delta/\delta_\Gamma}\big]$, the multiplicity of a degeneration graph~$\Gamma$ given by the degeneration formula is as follows:
\[ \int_{\prod_V\fvir{\M_V}{\delta_\Gamma}}
\ev^*\Delta\cup\prod_1^{n+g-1}\ev_i^*(\pt). \]
This integral splits as a product over the vertices of the graph once we use the K\"unneth decomposition of the diagonal. If $(b_i)$ is a basis and $\bigl(b_i^\#\bigr)$ the Poincar\'e dual basis, we have
\[
\Delta=\sum_i (-1)^{\mathrm{rk} b_i}b_i\otimes b_i^{\#}.
\] In our case, taking the basis $(1,\alpha,\beta,\pt)$ of $H^\bullet(E,\QQ)$, we have
\[ \Delta = 1\otimes\pt +\pt\otimes 1 + \alpha\otimes\beta -\beta\otimes\alpha. \]
Expanding, it amounts to sum over all the possible insertions of Poincar\'e dual classes of $E$ at the extremities of bounded edges of $\Gamma$, which we call a \textit{Poincar\'e insertion}. We now consider a~degeneration graph $\Gamma$ and aim at transforming it into a floor diagram. This is the content of Lemmas \ref{lem-type-vertices} and \ref{lem-complement-tree}.

\begin{lem}\label{lem-type-vertices}
 If $\Gamma$ is a degeneration graph with non-zero multiplicity, we have the following:
 \begin{enumerate}\itemsep=0pt
 \item[$(1)$] Inner vertices have genus $0$ or $1$, and a genus $1$ vertex is adjacent to a marked point.
 \item[$(2)$] A genus $0$ vertices is bivalent and its class is of the form $w\big[\PP^1\big]$.
 \item[$(3)$] For a marked genus $0$ vertex, adjacent Poincar\'e insertions are $1\in H^0(E,\QQ)$.
 \end{enumerate}
\end{lem}

\begin{proof}
 (1) Let $V$ be a vertex with genus $g_V$, valency $n_V$, and $m_V=0,1$ adjacent marked points. The dimension of $\fvir{\M_V}{\delta_\Gamma}$ is $n_V+m_V+g_V-1$. For the multiplicity associated to a Poincar\'e insertion to be non-zero, the sum of ranks of the insertions (Poincar\'e insertions coming from $\Delta$ and the $\pt$ if $m_V=1$) needs to match this dimension. The sum of ranks of these insertions is at most $(n_V-1)+2m_V$. The maximal rank of Poincar\'e insertions is $n_V-1$ and not $n_V$ due to the existence of a relation between the position of points at infinity: choosing only $\pt$ insertions at every flag yields $0$ invariant. Therefore, we have the inequality
 \[ n_V+m_V+g_V-1\leqslant n_V-1+2m_V. \]
 In other words, $g_V\leqslant m_V$, which implies (1): $g_V$ may only take the values $0,1$, and a genus $1$ vertex is marked.

 (2) Assume that $g_V=0$. As there all maps from a rational curve to an elliptic curve are constant, the class associated to $V$ is necessarily of the form $w\big[\PP^1\big]$. Fibers varying in a $1$-dimensional family, we may only impose a unique point constraint:
 \begin{itemize}\itemsep=0pt
 \item If $m_V=1$, we already have an interior point constraint, and the dimension is $(n_V+1)-1$. Thus, we have $n_V\leqslant 2$.
 \item If $m_V=0$, the dimension is $n_V-1$, and we the sum of rank of the Poincar\'e insertions is at most $1$, so that we also have $n_V\leqslant 2$.
 \end{itemize}
 In either case, as $n_V\geqslant 2$ due to the class intersecting both divisors of the singular locus, we get the bivalency statement from (2).

 (3) Assume that $V$ has genus $0$ and is marked ($m_V=1$). Then the dimension of $\M_V$ is $2$, and is already matched by the rank of the point insertion at the marked point, finishing the proof.
\end{proof}

Lemma \ref{lem-type-vertices} is a big step toward the transfiguration of a degeneration graph into a floor diagram: edges are already weighted, vertices are split into flat bivalent vertices and floors as advertised after Example \ref{expl-floor-diagram}. The total order comes from the total order on the components of~$\Y_0$. However, we still need three features:
 \begin{itemize}\itemsep=0pt
 \item delete unmarked genus $0$ vertices, which is achieved through Lemma \ref{lem-local-computation-bivalent-vertex},
 \item prove that the complement of genus $0$ marked vertices is a forest with each component containing a unique infinite vertex,
 \item the multiplicity matches $\bfm_\delta$.
 \end{itemize}
Consider $\Gamma$ our degeneration graph. We can cut $\Gamma$ at each genus $0$ marked vertices, replacing each bivalent vertex by two univalent vertices. The new graph (which may be disconnected) is denoted by $\widehat{\Gamma}$. We then remove the infinite vertices and get a non-compact graph $\widehat{\Gamma}^\circ$ called \textit{open cut graph}.

\begin{lem}\label{lem-complement-tree}
 If $\Gamma$ is a degeneration graph with non-zero multiplicity, and $\widehat{\Gamma}^\circ$ the associated open cut graph, we have the following:
 \begin{enumerate}\itemsep=0pt
 \item[$(1)$] Every connected component of $\widehat{\Gamma}^\circ$ has Euler characteristic $0$.
 \item[$(2)$] There are no compact connected component of $\widehat{\Gamma}^\circ$: every component of $\widehat{\Gamma}$ contains an infinite vertex.
 \item[$(3)$] $\widehat{\Gamma}$ has no cycle, and each connected component contains a unique infinite vertex.
 \item[$(4)$] The Poincar\'e insertions are uniquely determined:
 \begin{itemize}\itemsep=0pt
 \item for a marked genus $0$ vertex, $1,1\in H^\bullet(E,\QQ)$,
 \item for an unmarked genus $0$ vertex, $1,\pt\in H^\bullet(E,\QQ)$,
 \item for a (marked) genus $1$ vertex $1,\pt,\dots,\pt\in H^\bullet(E,\QQ)$.
 \end{itemize}
 \end{enumerate}
\end{lem}

\begin{proof}
 Let $C$ be a connected component of the compact graph $\widehat{\Gamma}$. We denote by
 \begin{itemize}\itemsep=0pt
 \item $\V_1$ the set of genus $1$ (inner) vertices,
 \item $\V_0$ the set of genus $0$ unmarked vertices,
 \item $\partial$ the set of univalent vertices resulting from the cut,
 \item $\partial_\infty$ the set of infinite vertices,
 \item $E_b$ the set of edges not adjacent to a vertex of $\partial\cup\partial_\infty$.
 \end{itemize}
 Notice that we counting the flags adjacent to the vertices of $\V_0\cup\V_1$ in two ways, we have that~${\sum_{\V_1}n_V +2|\V_0| = 2|E_b|+|\partial|+|\partial_\infty|}$.
 The dimension of $\prod_{V\in \V_0(C)\cup \V_1(C)}\fvir{\M_V}{\delta_\Gamma}$ is
 \begin{align*}
 \prod_{V\in \V_0(C)\cup \V_1(C)}\fvir{\M_V}{\delta_\Gamma} &=  \sum_{\V_1(C)} (n_V+1) + \sum_{\V_0(C)} 1  \\
& =   2|E_b|+|\partial|+|\partial_\infty| +|\V_1|-|\V_0| .
 \end{align*}
 For the invariant to be non-zero, it needs to match the sum of ranks of the insertions:
 \begin{itemize}\itemsep=0pt
 \item the points for each vertex of $\partial$ through Lemma \ref{lem-type-vertices},
 \item the point insertions for every vertex of $\V_1$,
 \item the diagonal insertions for every edge of $E_b$.
 \end{itemize}
 Therefore, we have
 \[ 2|E_b|+|\partial|+|\partial_\infty| +|\V_1|-|\V_0|= |\partial|+2|\V_1|+|E_b|. \]
 This yields
$|\V_1|+|\V_0|-|E_b| -|\partial_\infty|=0$.
 This quantity is precisely the Euler characteristic of $C$ where we deleted the vertices of $\partial_\infty$. Hence, we have (1).

 As the Euler characteristic is $0$, we have two possibilities:
 \begin{itemize}\itemsep=0pt
 \item the connected component is compact with a unique cycle,
 \item the connected component has no cycle and contains a unique end, meaning $|\partial_\infty|=1$.
 \end{itemize}

 The point (2) forbids the first possibility so that we have (3). To prove (2), we may use~\cite[Proposition 4.17]{blomme2024bielliptic} (and more precisely Step 5) that states that in the event where we have a~cycle, the invariant is $0$ since we have no monodromy in this situation.

 The statement about insertions follows from an induction on the trees of the forest $\widehat{\Gamma}^\circ$, pruning the branches until we are left with the edges adjacent to $\partial_\infty$.
\end{proof}

The local computation at the genus $1$ vertices has already been carried out in Theorem~\ref{theo-expression-local-correlated-invariant}. We end with the local computation at the genus $0$ vertices justifying the deletion of unmarked genus $0$ vertices before proving Theorem~\ref{theo-correspondence-diagrams}.

\begin{lem}\label{lem-local-computation-bivalent-vertex}
 For bivalent vertices, we have the following invariants:
 \[ \langle\langle\pt_0,1_w,1_{-w}\rangle\rangle^w_{0,0,(w,-w)} = 1\cdot (0)
 \qquad \text{and}\qquad
\langle\langle\pt_w,1_{-w}\rangle\rangle^w_{0,0,(w,-w)} = \frac{1}{w}\cdot (0). \]
\end{lem}

\begin{proof}
 Curves in the class $w\big[\PP^1\big]$ are mapped to some fiber of $E\times\PP^1$. Therefore, the only correlated invariant which is non-zero in this situation is for $(0)$. The uncorrelated computation has already been carried out in \cite[Lemma 3.3]{bousseau2021floor}: in both situations there is a unique map to the fiber, but with a $\ZZ/w\ZZ$ automorphism group in the second case.
\end{proof}

\begin{proof}[Proof of Theorem~\ref{theo-correspondence-diagrams}]
 The number of marked points is equal to $n+g-1$. We consider a degeneration $\Y_t$ of $E\times\PP^1$ with a central fiber $\Y_0$ having $n+g-1$ irreducible components. We apply the degeneration formula from Theorem~\ref{theo-refined-decomposition-surjective-case}, putting exactly one point constraint per component of the central fiber. Degeneration formula asserts that the correlated invariant we care about is a sum over the degeneration graphs. We need to show that the degeneration graphs with non-zero multiplicity are actually the floor diagrams involved in the sum, and that they share the same multiplicity.

 Given a degeneration graph $\Gamma$, Lemma \ref{lem-type-vertices} ensures that the vertices have the right genus and valency while Lemma \ref{lem-complement-tree}\,(3) states that the complement of all flat vertices has no cycle and each component has a unique end, so that $\Gamma$ is indeed a floor diagram. In particular, there are no automorphisms.

 Furthermore, by Lemma \ref{lem-complement-tree}\,(4), the Poincar\'e insertions coming from the K\"unneth decomposition of the diagonal are uniquely determined, so that the multiplicity splits as a product over the vertices. The multiplicity at genus $0$ vertices is provided by Lemma \ref{lem-local-computation-bivalent-vertex} and tells us that we may forget about genus $0$ unmarked vertices, as their multiplicity cancels with a factor of $\prod_e w_e$ when merging the adjacent edges.

 The local multiplicity at the genus $1$ vertices is provided by Theorem~\ref{theo-expression-local-correlated-invariant}, yielding the multiplicity $\bfm_\delta$.
\end{proof}

\begin{remarkk}
 The diagram setting can be adapted to deal with invariants including boundary constraints. To simplify notations, we only consider the case of interior point constraints.
\end{remarkk}

 \subsection{Quasi-modularity}
 \label{sec-modularity}

We now use the floor diagrams to prove some regularity statements for the invariants, starting with the generating series in the base direction: varying $a$ in the curve class $a[E]+b\big[\PP^1\big]$. For that purpose, we assume that $L=\O$, so that we may choose $\lambda_0=0$ and the term $(a\lambda_0)$ in the diagram multiplicities disappears.

We consider quasi-modular forms with values in the vector space $\QQ[\Tor_\delta(E)]$. Usual quasi-modular definitions extend naturally to the setting of functions in a vector space.

\begin{defi}
 Let $V$ be a complex vector space and $f\colon \HH\to V$ a meromorphic function on the Poincar\'e half-plane and $\Gamma\subset {\rm SL}_2(\ZZ)$ a subgroup.
 \begin{enumerate}\itemsep=0pt
 \item[(1)] We say that $f$ is modular of weight $k$ if for any $\abcd\in \Gamma$ we have
 \[ (c\tau+d)^kf\left(\frac{a\tau+b}{c\tau+d}\right) = f(\tau). \]
 \item[(2)] We say that $f$ is quasi-modular of weight $k$ and depth $s$ if there exists functions $f_j\colon \HH\to V$ such that for any $\abcd\in \Gamma$ we have \[ (c\tau+d)^k f\left(\frac{a\tau+b}{c\tau+d}\right) = \sum_0^s f_j(\tau)\left(\frac{c}{c\tau+d}\right)^j. \]
 \end{enumerate}
\end{defi}

A function $f\colon \HH\to V$ is (quasi-)modular if and only if all its coordinate functions are \mbox{(quasi-)modular}. Furthermore, in the quasi-modularity case, all $f_j$ are also quasi-modular forms and $f_0=f$. Furthermore, if the vector space $V$ is endowed with an $\CC$-algebra structure, the product of quasi-modular forms is also quasi-modular. If $\left(\begin{smallmatrix} 1 & 1 \\ 0 & 1 \\ \end{smallmatrix}\right)\in\Gamma$, we can set $\sfq={\rm e}^{2{\rm i}\pi\tau}$ and develop the quasi-modular forms in Fourier series, or conversely consider generating series in $\sfq$ as meromorphic functions potentially quasi-modular.

\begin{expl}
 We already saw that the generating series $\sum_{a=1}^\infty \boldsymbol{\sigma}_\delta(a)\sfq^a$ are quasi-modular forms for the congruence subgroup $\Gamma_0(\delta)$. Its derivatives $\sum_{a=1}^\infty a^n\boldsymbol{\sigma}_\delta(a)\sfq^a$ are as well.
\end{expl}

\begin{theo}\label{theo-quasi-modularity}
Let $\bfw$ be a tangency profile of length $n$. The generating series
\[ \sum_{a=1}^\infty \gen{\gen{\pt^{n+g-1}}}_{g,a[E],\bfw}^\delta \sfq^a \]
is a quasi-modular form for the congruence subgroup $\Gamma_0(\delta)$ with values in the group algebra $\CC[\Tor_\delta(E)]$.
\end{theo}

\begin{proof}
Up to the labeling of floors by $(a_V)$, there is a finite number of floor diagrams. It thus suffices to prove the statement for each unlabelled floor diagram.

Let $\Dfk$ a floor diagram without floor labels and let \smash{$W=\prod_{E_b(\Dfk)}w_e\prod_{E_\circ(\Dfk)} w_e^2$}. Let $\Dfk\bigl((a_V)\bigr)$ be the floor diagram obtained by labeling the vertices with $(a_V)$. As $\dfk\big[\frac{1}{\delta/\delta_\Gamma}\big]$ is compatible with the product, we have the following factorization
\begin{align*}
 \sum_{a=1}^\infty \biggl(\sum_{\Sigma a_V=a} \bfm_\delta(\Dfk((a_V)))\biggr) \sfq^a ={} & \sum_{(a_V)}\boldsymbol{m}_\delta(\Dfk\bigl((a_V)\bigr))\sfq^{\sum a_V} \\
 = {}& W\sum_{(a_V)}\d{\frac{1}{\delta/\delta_\Dfk}}\biggl(\prod_V a_V^{n_V-1}\boldsymbol{\sigma}_{\delta_\Dfk}(a_V)\sfq^{a_V}\biggr) \\
 = {}& W\cdot \d{\frac{1}{\delta/\delta_\Dfk}}\prod_V\left(\sum_{a_V=1}^\infty a_V^{n_V-1}\boldsymbol{\sigma}_{\delta_\Dfk}(a_V) \sfq^{a_V}\right).
\end{align*}
As each series is a quasi-modular form, so is their product.
\end{proof}

\begin{remarkk}
 The optimality of $\delta$ in the statement is ensured by the existence of floor diagrams with gcd $\delta$.
\end{remarkk}

\begin{remarkk}
 If we do not assume that $L=\O$ anymore, the coefficients of the generating series are multiplied by $(a\lambda_0)$, where $\delta\lambda_0=\lambda$. The naive version of quasi-modularity is not satisfied anymore. It is not really clear which kind of regularity to expect for these functions, with values in $\CC[E]$ and not the finite-dimensional $\CC[\Tor_\delta(E)]$.
\end{remarkk}

\subsection{Piecewise polynomiality}

We now vary the tangency orders, considering them as variables. We study the function
\[ N^\delta_{a,g}\colon\ \bfw = (w_1,\dots,w_n)\longmapsto \langle\langle\pt^{n+g-1},1_{w_1},\dots,1_{w_n}\rangle\rangle^\delta_{g,a[E],\bfw} \in \QQ[\Tor_\delta(E)]. \]
For this definition to make sense, we need to restrict to the sublattice where all the $w_i$ are divisible by $\delta$, so that the $\delta$-refinement makes sense.

Up to the labeling of edges by their weight, there is a finite number of floor diagrams of genus~$g$ with $n$ ends and with a class of the form $a[E]+\ast \big[\PP^1\big]$. The orientation is still part of the data, only the edge weights are missing. Let $\Dfk$ be one of them. To make $\Dfk$ into a~true floor diagram, one only needs to add weights to the edges. We denote by $\Omega_\Dfk$ the set of weightings~$\omega$ of~$\Dfk$, i.e., functions from edges to positive integers that make the floors and flat vertices balanced. For a~weighting $\omega\in\Omega_\Dfk$, we denote by $\delta_\omega$ the gcd between $\delta$ and its coordinates. The tangency profile is the restriction of the weighting to the infinite ends. For a fixed tangency profile $\bfw$, we denote by $\Omega_\Dfk(\bfw)$ the set of weightings that induce $\bfw$ on the ends of $\Dfk$. Although by assumption~$\bfw$ is divisible by $\delta$, $\omega$ may not be.

The multiplicity from Definition \ref{defi-multiplicity-diagram} contains as a factor a monomial in the edge weights, which we denote by $f_\Dfk$,
\[ f_\Dfk \colon\ \omega\in\Omega_\Dfk\longmapsto \prod_{E_b(\Dfk)}w_e\prod_{E_\circ(\Dfk)} w_e, \]
where the first product is over bounded edges of $\Dfk$ and the second product is over edges not adjacent to a flat vertex. The remaining part of the multiplicity depends on the floors weights~$a_V$ but also on the weighting $\omega$ through its gcd $\delta_\omega$. We denote it by
\[ \Phi_\Dfk(\delta_\omega)=\d{\frac{1}{\delta/\delta_\omega}}\bigg(\prod_V a_V^{n_V-1}\boldsymbol{\sigma}_{\delta_\omega}(a_V)\bigg), \]
so that $\boldsymbol{m}_\delta\bigl(\Dfk(\omega)\bigr)=\Phi_\Dfk(\delta_\omega)f_\Dfk(\omega)$. Even if we only care about the weighting that make the vertices balanced, this function is defined for all assignations of weights to the edges, not only the elements of $\Omega_\Dfk$. Before going to the main statement of the section, we reformulate the expression of $\Phi_\Dfk(\delta_\omega)$ to suit our purpose.

\begin{lem}\label{lem-multiplicity-polynomiality-statement}
 For a chosen diagram $\Dfk$ up to edge weights, there exists functions $\Upsilon^\Dfk_d$ such that
 \[ \Phi_\Dfk(\delta_\omega) = \sum_{d|\delta_\omega} \Upsilon^\Dfk_d(\delta_\omega)\vartheta_{\delta/d}. \]
\end{lem}

\begin{proof}
 We use the expression of $\boldsymbol{\sigma}_{\delta_\omega}$ highlighting the role of the $\vartheta_d$
 \[ \boldsymbol{\sigma}_{\delta_\omega}(a_V) = \sum_{d|\delta_\omega}\Upsilon_d^{\delta_\omega}(a_V)\vartheta_{\delta_\omega/d}. \]
 As the map $\dfk\big[\frac{1}{\delta/\delta_\omega}\big]$ is compatible with the product and $\dfk\big[\frac{1}{\delta/\delta_\omega}\big](\vartheta_{\delta_\omega/d})=\vartheta_{\delta/d}$, we have that
 \[ \Phi_\Dfk(\delta_\omega) = \prod_V \biggl(a_V^{n_V-1}\sum_{d_V|\delta_\omega}\Upsilon_{d_V}^{\delta_\omega}(a_V)\vartheta_{\delta/d_V}\biggr). \]
 We expand the product and use that $\vartheta_{\delta/d_1}\vartheta_{\delta/d_2}=\vartheta_{\delta/\gcd(d_1,d_2)}$,
 \[ \Phi_\Dfk(\delta_\omega) = \left(\prod a_V^{n_V-1}\right)\sum_{d|\delta_\omega}\biggl(\sum_{\substack{d_V|\delta_\omega \\ \gcd(\delta_V)=d}}\prod_V\Upsilon_{d_V}^{\delta_\omega}(a_V)\biggr)\vartheta_{\delta/d}. \]
 We thus have the existence and an expression for the desired coefficients $\Upsilon^\Dfk_d(\delta_\omega)$, which only depend on $\delta_\omega$ and $d$ once $\Dfk$ is fixed.
\end{proof}

Our statement relies on \cite[Theorem 4.2]{ardila2017double} which is a generalization of \cite[Theorem 1]{sturmfels1995vector}. To apply the theorem, we provide a second description of the weighting set $\Omega_\Dfk(\bfw)$ as the set of lattice points of a flow polytope. Let $A_\Dfk$ be the adjacency matrix of $\Dfk$: rows are indexed by vertices (floors, flat vertices, sinks, sources) and columns by edges. Coefficients are given by the following rule:
\[ A_\Dfk(V,e) = \begin{cases}
 \hphantom{-}0 &\text{if }V\notin e,\\
 \hphantom{-}1 & \text{if }e\text{ ends at }V,\\
 -1 & \text{if }e\text{ starts at }V.
\end{cases} \]
In particular, for an element \smash{$\omega\in\ZZ^{n+|E_b(\Dfk)|}$}, the image $A_\Dfk\omega$ is exactly the divergence at each vertex: the difference between incoming and outgoing weights at $V$. Therefore, if $\mathbf{d}=(d_V)$ is the vector with coordinate $|w_i|$ for an infinite vertex and $0$ else, we have
\[ \Omega_\Dfk(\bfw) = \big\{ \omega\in\NN^{n+ E_b(\Dfk)} \text{ s.t. } A_\Dfk\omega=\mathbf{d} \big\}. \]

\begin{theo}\label{theo-quasi-polynomiality}
For fixed $a$, $g$, $n$, $\delta$, the function $\bfw\mapsto N^\delta_{a,g}(\bfw)$ is piecewise polynomial in the sense that there exists piecewise polynomial functions $P_d(\bfw)$ for $d|\delta$ such that
\[ N^\delta_{a,g}(\bfw) = \sum_{d|\delta}P_d(\bfw)\vartheta_d. \]
\end{theo}

\begin{proof}
Up to the weighting of the edges, there is a finite number of floor diagrams of genus $g$ with class of the form $a[E]+\ast \big[\PP^1\big]$ and $n$ ends. Grouping them together, we get
\[ N^\delta_{a,g}(\bfw) = \sum_\Dfk N^\delta_\Dfk(\bfw),\qquad \text{with}\quad N^\delta_\Dfk(\bfw) = \sum_{\omega\in\Omega_\Dfk(\bfw)} \Phi_\Dfk(\delta_\omega) f_\Dfk(\omega). \]
Thus, to prove piecewise polynomiality, it is sufficient to consider a unique unweighed floor diagram $\Dfk$.

Let $\Dfk$ be such an unweighed floor diagram. If there was no $\Phi_\Dfk(\delta_\omega)$, \cite[Theorem 4.2]{ardila2017double} already applies and yields piecewise polynomiality. Due to the presence of this factor depending on the gcd $\delta_\omega$, our second step is thus to first use the expression of $\Phi_\Dfk(\delta_\omega)$, split the sum and perform a \textit{multiplicative summation by parts} (as described before Theorem~\ref{theo-expression-local-correlated-invariant}). We aim to change the sum over $\omega$ and a given $\delta_\omega$ into a sum over the $\omega$ in a sublattice. Using Lemma \ref{lem-multiplicity-polynomiality-statement}, we have
\[ N^\delta_\Dfk(\bfw) = \sum_{\omega\in\Omega_\Dfk(\bfw)} \biggl(\sum_{d|\delta_\omega}\Upsilon^\Dfk_d(\delta_\omega)\vartheta_{\delta/d}\biggr) f_\Dfk(\omega). \]
To perform the multiplicative summation by parts, we can find coefficients $\boldsymbol{\gamma}^\Dfk_d$, which are linear combinations of the $\vartheta_{\delta/d}$, such that for each $\delta_\omega$ we have
\[ \sum_{d|\delta_\omega}\boldsymbol{\gamma}^\Dfk_d =
\sum_{d|\delta_\omega}\Upsilon^\Dfk_d(\delta_\omega)\vartheta_{\delta/d}. \]
Such coefficients exist because the equations given for each $\delta_\omega$ form a triangular linear system. We can now switch the two sums and get rid of the dependence in the gcd $\delta_\omega$
\[ N^\delta_\Dfk(\bfw) =
\sum_{\omega\in\Omega_\Dfk(\bfw)} \biggl(\sum_{d|\delta_\omega}\boldsymbol{\gamma}^\Dfk_d\biggr) f_\Dfk(\omega)
=\sum_{d|\delta} \boldsymbol{\gamma}^\Dfk_d \sum_{\substack{\omega\in\Omega_\Dfk(\bfw) \\ d|\delta_\omega}} f_\Dfk(\omega). \]

The last step is now to apply \cite[Theorem 4.2]{ardila2017double}. For each $d|\delta$, we can rewrite
\[ \sum_{\substack{\omega\in\Omega_\Dfk(\bfw) \\ d|\omega}} f_\Dfk(\omega) = \sum_{\omega'\in\Omega_\Dfk(\bfw/d)} f_\Dfk(d\omega'). \]
The function is piecewise quasi-polynomial with respect to the chamber complex of $A_\Dfk$. Furthermore, as $A_\Dfk$ is actually unimodular (it is a submatrix of the root system $A_{n-1}$, see \cite[Example~4.4]{ardila2017double}), we in fact get piecewise polynomiality and not only quasi-polynomiality. Summing over all diagrams, we conclude.
\end{proof}

\begin{expl}
 Assume $\delta=p$ is a prime number and $\Dfk$ is a diagram up to edge weights. Then we have two kinds of weightings:
 \begin{itemize}\itemsep=0pt
 \item the weightings with gcd $1$ and multiplicity $\Upsilon^\Dfk_1(1)\vartheta_p \cdot f_\Dfk(\omega)$,
 \item the weightings with gcd $p$ and multiplicity $\left(\Upsilon^\Dfk_1(p)\vartheta_p + \Upsilon^\Dfk_p(p)\vartheta_1\right)\cdot f_\Dfk(\omega)$.
 \end{itemize}
 The multiplicative summation by parts tells us to count all weightings with the multiplicity of gcd $1$ weightings, and correct by counting the gcd $p$ weightings with multiplicity
 \[ \bigl(\bigl(\Upsilon^\Dfk_1(p)-\Upsilon^\Dfk_1(1)\bigr)\vartheta_p + \Upsilon^\Dfk_p(p)\vartheta_1\bigr)f_\Dfk(\omega). \]
\end{expl}

\begin{figure}[th]
 \centering
\includegraphics{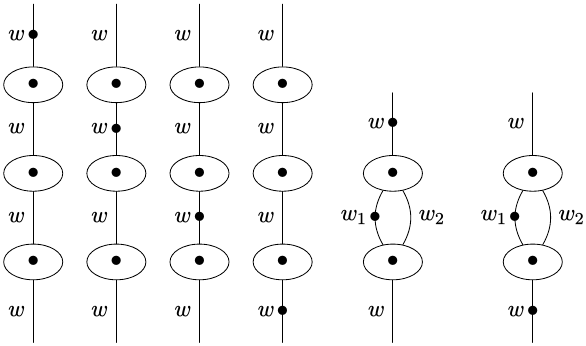}

 \caption{Floor diagrams of genus $3$.}
 \label{fig-genus-3-diagrams}
\end{figure}

\begin{expl}
 We end this section with an explicit computation. We take $n=2$ and $g=3$ and a chosen $\delta$, for instance $\delta=2$. As $n=2$, the tangency profile is necessarily $(w,-w)$ and~$w$ is divisible by $\delta$. Up to weights of the edges and floors, we only have six floor diagrams, of two different kinds. For the first four floor diagrams, the multiplicity is a constant term in~$\QQ[\Tor_\delta(E)]$ times $w^9$, which is indeed a polynomial in $w$. We now care about the diagrams of the second kind. To choose a weighting of the edges, one needs to choose a splitting $w=w_1+w_2$. If they are both even, the gcd $\delta_\omega$ of the diagram is $2$, otherwise both of them are odd and the gcd is $1$. The multiplicity is
 \[ a_1^2a_2^2\d{\frac{1}{2/\delta_\omega}}\left(\boldsymbol{\sigma}_{\delta_\omega}(a_1)\boldsymbol{\sigma}_{\delta_\omega}(a_2)\right)w^3w_1^2w_2^3, \]
 where we have the two expressions
 \begin{align*}
 \boldsymbol{\sigma}_1(a) = \sigma(a)\vartheta_1, \qquad
 \boldsymbol{\sigma}_2(a) = \overline{\sigma}^2(a)\vartheta_1 + (\sigma(a)-\overline{\sigma}^2(a))\vartheta_2.
 \end{align*}
 Therefore, in the odd case the multiplicity is
$a_1^2a_2^2\sigma(a_1)\sigma(a_2)\vartheta_2\cdot w^3w_1^2w_2^3$, \
 and in the even case
 \[ a_1^2a_2^2\bigl(\overline{\sigma}^2(a_1)\overline{\sigma}^2(a_2)\vartheta_1 + \bigl(\sigma(a_1)\sigma(a_2) - \overline{\sigma}^2(a_1)\overline{\sigma}^2(a_2)\bigr)\vartheta_2 \bigr)\cdot w^3w_1^2w_2^3. \]
 To get the result, we sum the multiplicity in the odd case over $1\leqslant w_1\leqslant w-1$ to get a first polynomial, and sum the difference between the two multiplicity over the even $w_1$ satisfying the same inequalities.
\end{expl}

\subsection*{Acknowledgements}
 T.B.\ is supported by the SNF grant 204125; F.C.\ is supported by the Ambizione grant PZ00P2-208699/1. F.C.\ is also partially supported by the MIUR Excellence Department Project MatMod@TOV, CUP E83C23000330006, awarded to the Department of Mathematics, University of Rome Tor Vergata, and also acknowledges the support of the PRIN Project ``Moduli spaces and birational geometry'' 2022L34E7W.
The authors would like to thank D.~Ranganathan, A.~Kumaran, S.~Molcho for helpful discussion around the subject of this paper.

\pdfbookmark[1]{References}{ref}
\LastPageEnding

\end{document}